%% file: 2005-11.tex
\documentclass{gtart_h}  

\input gtspec

\lognumber{538}
\received{7 December 2004}
\volumenumber{9}\papernumber{11}\volumeyear{2005}
\pagenumbers{341}{373}   
\published{1 March 2005}
\accepted{28 February 2005}
\proposed{Tomasz Mrowka}
\seconded{Ralph Cohen, Gang Tian}

\def\S{Section }
\usepackage{amsmath,amssymb,amscd,verbatim}
\usepackage[all]{xy}

\newtheorem{thm}{Theorem}[section]
\newtheorem{corollary*}{Corollary}
\newtheorem{lemma}{Lemma}
\newtheorem{proposition}{Proposition}

\theoremstyle{definition}

\newtheorem*{definition*}{Definition}
\newtheorem*{rem*}{Remark}
\newtheorem*{proof*}{Proof}
\newtheorem{prenl*}{Preliminaries}
\newtheorem{examples*}{Examples}

\newcommand\ev{\operatorname{ev}}
\newcommand\PU{\operatorname{PU}}
\newcommand\SU{\operatorname{SU}}
\newcommand\UG{\operatorname{U}}
\newcommand\proj{\operatorname{proj}}
\newcommand\Cstar{$C^*\kern-2pt$}
\newcommand\Id{\operatorname{Id}}

\newcommand\CC{\mathbb{C}}

\newcommand\fL{L^{-1}}

\newcommand\RR{\mathbb{R}}

\newcommand\ZZ{\mathbb{Z}}

\newcommand\cJ{\mathcal{J}}
\newcommand\Cinf{\mathcal{C}^\infty}

\newcommand\Diag{\operatorname{Diag}}

\newcommand\Hom{\operatorname{Hom}}
\newcommand\Ad{\operatorname{Ad}}

\newcommand\R{\mathcal{R}}

\newcommand{\Aut}{\operatorname{Aut}}
\newcommand{\Br}{\operatorname{Br}}
\newcommand{\C}{\mathbb{C}}
\newcommand{\E}{\widetilde{\mathcal E}}

\newcommand{\cH}{{\mathcal H}}
\newcommand{\cK}{{\mathcal K}}
\newcommand{\cA}{{\mathcal A}}

\newcommand{\Min}{\text{ in }}

\newcommand{\Fred}{\operatorname{Fred}}
\newcommand{\TC}{\Theta}

\newcommand{\TZ}{\theta}

\newcommand\bbQ{\mathbb{Q}}
\newcommand\bbZ{\mathbb{Z}}
\newcommand\bbN{\mathbb{N}}
\newcommand\bbC{\mathbb{C}}

\def\Kc{K_c}

\def\C{\mathbb C}
\def\R{\mathbb R}

\def\Z{\mathbb Z}
\def\Q{\mathbb Q}
\def\E{\mathcal E}
\def\K{\mathcal K}

\def\B{\mathcal B}

\def\Hom{\operatorname{Hom}}

\def\Tsn{\operatorname{tor}}

\def\index{\operatorname{index}}
\def\End{\operatorname{End}}
\newcommand\spinC{\text{spin${}^\bbC$}}
\def\Cl{\operatorname{Cl}}
\def\cli{\operatorname{cl}}

\def\GL{\operatorname{GL}}

\def\A{{\mathcal A}} \def\B{{\mathcal B}}

\def\H{\mathcal H}

\def\E{\mathcal E}
\def\U{{\mathcal U}}

\def\<{\langle}
\def\>{\rangle}

\newcommand{\Ch}{\operatorname{Ch}}

\newcommand{\nc}{\newcommand}
\nc{\nt}{\newtheorem}
\nc{\gf}[2]{\genfrac{}{}{0pt}{}{#1}{#2}}
\nc{\mb}[1]{{\mbox{$ #1 $}}}
\nc{\real}{{\mathbb R}}
\nc{\comp}{{\mathbb C}}
\nc{\ints}{{\mathbb Z}}
\nc{\Ltoo}{\mb{L^2({\mathbf H})}}
\nc{\rtoo}{\mb{{\mathbf R}^2}}
\nc{\slr}{{\mathbf {SL}}(2,\real)}
\nc{\slz}{{\mathbf {SL}}(2,\ints)}
\nc{\su}{{\mathbf {SU}}(1,1)}
\nc{\so}{{\mathbf {SO}}}
\nc{\hyp}{{\mathbb H}}
\nc{\disc}{{\mathbf D}}
\nc{\torus}{{\mathbb T}}

\nc{\ca}{{\mathcal A}}
\nc{\cag}{{{\mathcal A}^\Gamma}}
\nc{\cg}{{\mathcal G}}
\nc{\chh}{{\mathcal H}}
\nc{\ck}{{\mathcal B}}
\nc{\cm}{{\mathcal M}}
\nc{\cs}{{\mathcal S}}
\nc{\cz}{{\mathcal Z}}
\nc{\sind}{\sigma{\rm -ind}}

\begin{document}
\abovedisplayskip=6pt plus3pt minus6pt
\belowdisplayskip=6pt plus3pt minus6pt
\parskip=6pt plus3pt minus3pt

\title{The index of projective families of elliptic operators}
\authors{Varghese Mathai\\Richard B Melrose\\Isadore M Singer}
\shortauthors{Varghese Mathai, Richard B Melrose\\Isadore M Singer}
\addresses{Department of Pure Mathematics, University of Adelaide\\
Adelaide 5005, Australia\\
\smallskip\\
Department of Mathematics, Massachusetts Institute of Technology\\
Cambridge, Mass 02139, USA\\
\smallskip\\
Department of Mathematics, Massachusetts Institute of Technology\\
Cambridge, Mass 02139, USA\\
\medskip\\
{\rm Email:} {\tt\mailto{vmathai@maths.adelaide.edu.au}, 
\mailto{rbm@math.mit.edu}, \mailto{ims@math.mit.edu}}}
\asciiemail{vmathai@maths.adelaide.edu.au, rbm@math.mit.edu,
ims@math.mit.edu}

\asciiaddress{Department of Pure Mathematics, University of Adelaide\\
Adelaide 5005, Australia\\
Department of Mathematics, Massachusetts Institute of Technology\\
Cambridge, Mass 02139, USA\\
Department of Mathematics, Massachusetts Institute of Technology\\
Cambridge, Mass 02139, USA}

\asciikeywords{Projective vector bundles, twisted K-theory, projective families
of elliptic operators, Index theorem, determinant lines, twisted Chern character}

\begin{abstract} An index theory for projective families of elliptic
pseudodifferential operators is developed. The topological and the analytic
index of such a family both take values in twisted $K$--theory of the
parametrizing space, $X.$ The main result is the equality of these two notions
of index when the twisting class is in the torsion subgroup of
$H^3(X;\bbZ).$ The Chern character of the index class is then
computed.
\end{abstract}

\asciiabstract{An index theory for projective families of elliptic
pseudodifferential operators is developed. The topological and the analytic
index of such a family both take values in twisted K-theory of the
parametrizing space, X.  The main result is the equality of these two notions
of index when the twisting class is in the torsion subgroup of
H^3(X;Z).  The Chern character of the index class is then
computed.}

\maketitle

\section*{Introduction}\addcontentsline{toc}{section}{Introduction}

In this paper we develop an index theory for projective families of
elliptic pseudodifferential operators. Such a family, $\{D_b, \;b\in X\},$
on the fibers of a fibration \begin{equation} \phi \co M\longrightarrow X
\label{mms2001.320}\end{equation} with base $X,$ and typical fibre $F,$ is
a collection of local elliptic families, for an open covering of the base,
acting on finite-dimensional vector bundles of fixed rank where the usual
compatibility condition on triple overlaps, to give a global family, may
fail by a scalar factor. These factors define an integral 3--cohomology
class on the base, the Dixmier--Douady class, $\TC\in H^3(X,\bbZ).$ We show
that both the analytic and topological index of $D_b$ may be defined as
elements of twisted $K$--theory, with twisting class $\TC,$ and that they
are equal. In this setting of finite-dimensional bundles the twisting class
is necessarily a \emph{torsion} class. We also compute the Chern character
of the index in terms of characteristic classes. When the twisting class
$\TC$ is trivial, these results reduce to the Atiyah--Singer index theorem
for families of elliptic operators, \cite{MR43:5554}.

The vector bundles on which $D_b$ acts, with this weakened compatibility
condition, form a \emph{projective vector bundle}\footnote{Generally these
are called \emph{gauge bundles} in the physics literature}. In the torsion
case, elements of twisted $K$--theory may be represented by differences of
such projective vector bundles and, after stabilization, the local index
bundles of the family give such a difference and so define the analytic
index. Nistor and Troitsky have recently given a similar definition in
\cite{Nistor-Troisky1}. The topological index is defined, as in the
untwisted case, by push-forward from a twisted $K$ class of compact support
on the cotangent bundle. The proof we give of the equality of the analytic
and the topological index is a generalization of the axiomatic proof of the
families index theorem of \cite{MR43:5554}. This is done by proving that
all normalized, functorial index maps, in twisted $K$--theory, satisfying
excision and multiplicativity coincide and then showing that the analytic
index and the topological index satisfy these conditions.

Twisted $K$--theory arises naturally when one considers the Thom isomorphism
for a real Riemannian vector bundle, $E,$ in the even-rank case and when
$X$ is compact. The compactly supported $K$--theory of $E,$ $\Kc(E),$ is
then isomorphic to the twisted $K$--theory of $X,$ $K(C(X,\Cl(E)))$, where
$\Cl(E)$ denotes the Clifford algebra bundle of $E,$ and $C(X,\Cl(E))$
denotes the algebra of continuous sections. There is a similar statement
with a shift in degree for odd-dimensional vector bundles. The twisting of
families of Dirac operators by projective vector bundles provides many
examples of elliptic families and, as in the untwisted case, the Chern
character has a more explicit formula.

Recently, there has been considerable interest in twisted $K$--theory by
physicists, with elements of twisted $K$--theory interpreted as charges of
$D$--branes in the presence of a background field; cf \cite{MR87j:81081},
\cite{Witten3}, \cite{MR2000a:81190}.

There is an alternative approach to the twisted index theorem, not carried
out in detail here. This is to realize the projective family as an ordinary
equivariant family of elliptic pseudodifferential operators on an
associated principal $\PU(n)$--bundle. Then the analytic index and the
topological index are elements in the ${\rm U}(n)$--equivariant
$K$--theory of this
bundle. Their equality follows from the equivariant index theorem for
families of elliptic operators as in \cite{Tsuboi1}. The proof can then be
completed by showing that the various definitions of the analytic index and
of the topological index agree for projective families of elliptic
operators.

A subsequent paper will deal with the general case of the twisted index
theorem when the twisting 3--cocycle is not necessarily torsion. Then there
is no known finite-dimensional description of twisted $K$--theory, and even
to formulate the index theorem requires a somewhat different approach.

The paper is organized as follows. A review of twisted $K$--theory, with an
emphasis on the interpretation of elements of twisted $K$--theory as
differences of projective vector bundles is given in \S\ref{Review}. The
definition of general projective families of pseudodifferential operators
is explained in \S\ref{Analytic}, leading to the definition of the
analytic index, in the elliptic case, as an element in twisted
$K$--theory. The definition of the topological index is given in
\S\ref{Topological} and \S\ref{Equality} contains the proof of the
equality of these two indices. In \S\ref{Chern} the Chern character of the
analytic index is computed and in \S\ref{Determinant} the determinant
bundle is discussed in this context. Finally, \S\ref{Dirac} contains a brief
description of Dirac operators.

{\bf Acknowledgements}\qua VM acknowledges support from the 
Australian Research Council and the Clay Mathematics Institute,
RBM from NSF grant DMS-0104116 and IMS from DOE contract
DE-FG02-88ER25066.

\section{Review of twisted $K$--theory}\label{Review}

General references for most of the material summarized here are
\cite{Rosenberg2, Rosenberg1}.

\subsection{Brauer groups and the Dixmier--Douady invariant}

We begin by reviewing some results due to Dixmier and Douady,
\cite{Dixmier-Douady1}. Let $X$ be a smooth manifold, let $\cH$ denote an
infinite-dimensional, separable, Hilbert space and let $\cK$ be the
\Cstar--algebra of compact operators on $\cH.$ Let $\UG(\cH)$ denote the group
of unitary operators on $\cH$ endowed with the strong operator topology
and let $\PU(\cH)= \UG(\cH)/\UG(1)$ be the projective unitary group with the
quotient space topology, where $\UG(1)$ consists of scalar multiples of the
identity operator on $\cH$ of norm equal to $1.$ Recall that, if $G$ is a
topological group, principal $G$ bundles over $X$ are classified up to
isomorphism by the first cohomology of $X$ with coefficients in the sheaf
of germs of continuous functions from $X$ to $G,$ $H^1(X,\underline{G}),$
where the transition maps of any trivialization of a principal $G$ bundle
over $X$ define a cocycle in $Z^1(X, \underline{G})$ with fixed cohomology
class. The exact sequence of sheaves of groups,
$$
1\longrightarrow \underline{ \UG(1)}\longrightarrow
\underline{\UG(\cH)}\longrightarrow \underline{\PU(\cH)}\longrightarrow1
$$
gives rise to the long exact sequence of cohomology groups,
\begin{equation}\label{lexact1}
\dots\longrightarrow H^1(X, \underline{\UG(\cH)})\longrightarrow H^1(X,
\underline{\PU(\cH)}) \stackrel{\delta_1}{\longrightarrow } H^2(X,
\underline{\UG(1)})\longrightarrow 1 .
\end{equation}

Since $\UG(\cH)$ is contractible in the strong operator topology, the sheaf
$\underline{\UG(\cH)}$ is soft and the sheaf cohomology vanishes,
$H^1(X, \underline{\UG(\cH)}) =\{0\}.$ Equivalently, every Hilbert
bundle over $X$ is trivializable in the strong operator topology. In fact,
Kuiper \cite{Kuiper1} proves the stronger result that $\UG(\cH)$ is
contractible in the operator norm topology, so the same conclusion holds in
this sense too. It follows from \eqref{lexact1} that $\delta_1$ is an
isomorphism.  That is, principal $\PU(\cH)$ bundles over $X$ are classified
up to isomorphism by $H^2(X, \underline{\UG(1)}).$ From the exact sequence
of groups,
$$
1\longrightarrow \Z \longrightarrow \R\longrightarrow \UG(1)\longrightarrow 1
$$
we obtain the long exact sequence of cohomology groups,
$$
\cdots\longrightarrow H^2(X, \underline{\R}) \longrightarrow H^2(X,
\underline{\UG(1)}) \stackrel{\delta_2}{\longrightarrow} H^3(X, \Z)
\longrightarrow H^3(X,\underline{\R})\longrightarrow\cdots.
$$
Now $H^j(X,\underline{\R}) = \{0\}$ for $j>0$ since
$\underline{\R}$ is a fine sheaf, therefore $\delta$ is also
an isomorphism. That is, principal $\PU(\cH)$ bundles over $X$ are
also classified up to isomorphism by $H^3(X, \Z).$ The class
$\delta_2\delta_1([P])\in H^3(X,
\Z)$ is called the \emph{Dixmier--Douady class} of the principal
$\PU(\cH)$ bundle $P$ over $X$, where $[P] \in H^1(X,
\underline{\PU(\cH)}).$

For $g\in \UG(\cH),$ let $Ad(g)$ denote the automorphism $T\longrightarrow
gTg^{-1}$ of $\cK.$ As is well-known, $\Ad$ is a continuous homomorphism of
$\UG(\cH),$ given the strong operator topology, onto $\Aut(\cK)$ with
kernel the circle of scalar multiples of the identity where $\Aut(\cK)$ is
given the point-norm topology, that is the topology of pointwise
convergence of functions on $\K$, cf \cite{MR2000c:46108},
Chapter~1. Under this homomorphism we may identify $\PU(\cH)$ with
$\Aut(\cK).$ Thus

\begin{proposition}[Dixmier--Douady \cite{Dixmier-Douady1}]\label{mms2001.342} 
The isomorphism classes of locally trivial bundles over $X$ with fibre
$\cK$ and structure group $\Aut(\cK)$ are paramet\-rized by $H^3(X, \ZZ).$
\end{proposition}

Since $\cK\otimes\cK\cong \cK,$ the isomorphism classes of locally trivial
bundles over $X$ with fibre $\cK$ and structure group $\Aut(\cK)$ form a
group under the tensor product, where the inverse of such a bundle is the
conjugate bundle. This group is known as the \emph{infinite Brauer group}
and is denoted by $\Br^\infty(X)$ (cf \cite{Parker1}). So, essentially as
a restatement of Proposition~\ref{mms2001.342}
\begin{equation}
\Br^\infty(X) \cong H^3(X, \ZZ)
\label{mms2001.343}\end{equation}
where the cohomology class in $H^3(X, \ZZ)$ associated to a locally trivial
bundle $\mathcal E$ over $X$ with fibre $\cK$ and structure group
$\Aut(\cK)$ is again called the Dixmier--Douady invariant of $\mathcal E$ and
is denoted by $\delta(\mathcal E).$

In this paper, we will be concerned mainly with torsion classes in $H^3(X,
\Z).$ Let $\Tsn(H^3(X,\ZZ))$ denote the subgroup of torsion elements in
$H^3(X,\ZZ).$ Suppose now that $X$ is compact. Then there is a well-known
description of $\Tsn(H^3(X,\ZZ))$ in terms of locally trivial bundles of
finite-dimensional Azumaya algebras over $X,$
\cite{Donovan-Karoubi1}. Recall that an Azumaya algebra of rank $m$ is an
algebra which is isomorphic to the algebra of $m\times m$ matrices, 
$M_m(\bbC).$

\begin{definition*}\label{Azumaya}An \emph{Azumaya bundle} over a manifold
$X$ is a vector bundle with fibres which are Azumaya algebras and which has
local trivialization reducing these algebras to $M_m(\bbC).$
\end{definition*}

\begin{rem*}
The terminology `Azumaya bundle' is not the only one in the literature for
this notion. The oldest is that of ``$n$--homogeneous $C^*$--algebra''
due to Kaplansky, cf \cite{Kaplansky}, page 236.  These are
by definition $C^*$--algebras all of whose irreducible
representations have dimension equal to $n$. A highly nontrivial
theorem of Takesaki--Tomiyama, \cite{Takesaki-Tomiyama} 
shows that these algebras
are Azumaya algebras over a commutative algebra
and hence they identify with the space of continuous
sections of a bundle of finite dimensional simple algebras.
(A general reference is Diximier's book on $C^*$--algebras, \cite{Dixmier}. 
For the reader's benefit, we
recall that an Azumaya algebra $\cA$ over a ring $R$ is
an algebra that is a finitely generated, projective
algebra over $R$ such that $\cA \otimes \cA^{op} \simeq {\rm End}_R(\cA)$.)
\end{rem*}

An example of an Azumaya bundle over $X$ is the algebra, $\End(E),$ of all
endomorphisms of a vector bundle $E$ over $X.$ Two Azumaya bundles
$\mathcal E$ and $\mathcal F$ over $X$ are said to be equivalent if there
are vector bundles $E$ and $F$ over $X$ such that $\mathcal
E\otimes\End(E)$ is isomorphic to $\mathcal F\otimes\End(F).$ In
particular, an Azumaya bundle of the form $\End(E)$ is equivalent to $C(X)$
for any vector bundle $E$ over $X.$ The group of all equivalence classes of
Azumaya bundles over $X$ is called the Brauer group of $X$ and is denoted
by $\Br(X).$ We will denote by $\delta'(\mathcal E)$ the class in
$\Tsn(H^3(X,\ZZ))$ corresponding to the Azumaya bundle $\mathcal E$ over
$X.$ It is constructed using the same local description as above but with
$\mathcal{H}$ now finite dimensional. Serre's theorem, \cite{MR39:5586a}
asserts that
\begin{equation}
\Br(X)\cong\Tsn(H^3(X,\ZZ)).
\label{mms2001.344}\end{equation}
This is also the case if we consider smooth Azumaya bundles.

Another important class of Azumaya algebras arises from the bundles of
Clifford algebras of vector bundles. For a real vector bundle, with fibre
metric, the associated bundle of complexified Clifford algebras is an
Azumaya bundle for even rank and for odd rank is a direct sum of two
Azumaya bundles. In the even-rank case, the Dixmier--Douady invariant of
this Azumaya bundle is the third integral Stiefel--Whitney class,
the vanishing of which is equivalent to the existence of a $\spinC$
structure on the bundle.

Thus we see that there are two descriptions of $\Tsn(H^3(X, \ZZ)),$ one in
terms of Azumaya bundles over $X,$ and the
other as a special case of locally trivial bundles over $X$ with fibre
$\cK$ and structure group $\Aut(\cK).$ These two descriptions are related
as follows. Given an Azumaya bundle $\mathcal E$ over $X,$ the
tensor product $\mathcal E\otimes \cK$ is a locally trivial bundle over $X$
with fibre $M_m(\CC)\otimes\cK\cong \cK$ and structure group $\Aut(\cK),$
such that $\delta'(\mathcal E) = \delta(\mathcal E\otimes \cK)$. Notice
that the algebras $C(X, \mathcal E)$ and $C(X, \mathcal E\otimes \cK) =
C(X, \mathcal E)\otimes \cK$ are Morita equivalent. Moreover if $\mathcal
E$ and $\mathcal F$ are equivalent as Azumaya bundles over $X$
then $\mathcal E\otimes \cK$ and $\mathcal F\otimes \cK$ are isomorphic
over $X,$ as locally trivial bundles with fibre $\cK$ and structure group
$\Aut(\cK).$ To see this, recall that by the assumed equivalence there are
vector bundles $E$ and $F$ over $X$ such that $\mathcal E\otimes\End(E)$ is
isomorphic to $\mathcal F\otimes\End(F).$ Tensoring both bundles with
$\cK(\cH),$ we see that $\mathcal E\otimes {\cK}(E\otimes\cH)$ is
isomorphic to $\mathcal F\otimes {\cK}(F\otimes\cH)$, where
$\cK(E\otimes\cH)$ and $\cK(F\otimes\cH)$ are the bundles of compact
operators on the infinite dimensional Hilbert bundles $E\otimes\cH$ and
$F\otimes\cH$, respectively.  By the contractibility of the unitary group
of an infinite dimensional Hilbert space in the strong operator topology,
the infinite dimensional Hilbert bundles $E\otimes\cH$ and $F\otimes\cH$
are trivial, and therefore both $\cK(E\otimes\cH)$ and $\cK(F\otimes\cH)$
are isomorphic to the trivial bundle $X\times \cK.$ It follows that
$\mathcal E$ and $\mathcal F$ are equivalent Azumaya bundles
over $X,$ if and only if $\mathcal E\otimes\cK$ and $\mathcal F \otimes
\cK$ are isomorphic, as asserted.

Recall that a \Cstar--algebra $A$ is said to be \emph{stably unital} if
there is a sequence of projections $p_n\in A\otimes \K$ such that in the
strong topology $Tp_n\longrightarrow T$ for each $T\in A\otimes \K.$ In
particular $\K$ itself is stably unital since each compact operator can be
approximated by finite rank operators. It follows that any unital algebra
is stably unital.

\begin{lemma}\label{stably_unital} If $X$ is a compact manifold and
$\mathcal E$ is a locally trivial bundle over $X$ with fibre $\cK$ and
structure group $\Aut(\cK)$ then $C(X,\E),$ the \Cstar--algebra of
continuous sections of ${\mathcal E},$ is stably unital if and only if its
Dixmier--Douady invariant is a torsion element in $H^3(X,\Z).$
\end{lemma}

\begin{proof} The assumption that $C(X, \E)$ is stably unital where $\E$
has fibres isomorphic to $\K,$ implies in particular that there is a
non-trivial projection $p\in C(X, \E).$ This must be of finite rank in each
fibre. The \Cstar--algebra $pC(X, \E)p$ is a corner in $C(X, \E)$ in the
sense of Rieffel, \cite{MR84k:46045}, so these algebras are Morita
equivalent. They are both continuous trace \Cstar--algebras with the same
spectrum, which is equal to $X,$ and therefore the same Dixmier--Douady
invariant by the classification of continuous trace \Cstar--algebras
\cite{MR2000c:46108}. Since $pC(X,\E)p$ is the \Cstar--algebra of sections
of an Azumaya bundle over $X,$ the Dixmier--Douady invariant is a torsion
element in $H^3(X, \Z).$

The converse is really just Serre's theorem; given $\TC\in\Tsn(H^3(X, \Z))$
there is a principal $\PU(m)$ bundle over $X$ whose Dixmier--Douady
invariant is $\TC$ where the order of $\TC$ necessarily divides $m.$ The
associated Azumaya bundle $\cA$ also has the same
Dixmier--Douady invariant.  Then $\E = \cA\otimes \K$ is a locally trivial
bundle over $X$ with fibre $\cK$ and structure group $\Aut(\cK)$ and with
Dixmier--Douady invariant $\delta(\E) = \TC.$ So there is a non-trivial
projection $p_1 \in C(X, \E)$ such that $p_1 C(X, \E)p_1 = C(X, \cA).$ In
fact, one can define a nested sequence of Azumaya bundles, $\cA_j,$
over $X$ defined by $\cA_j = \cA\otimes M_j(\C)$ together with the
corresponding projections $p_j \in C(X, \E)$ such that $p_j C(X, \E)p_j =
C(X, \cA_j)$ for all $j\in \mathbb N.$ Then $\{p_j\}_{j\in \mathbb N}$ is
an approximate identity of projections in $C(X, \E)$, that is, $C(X, \E)$
is stably unital.
\end{proof}

\begin{rem*}\label{mms2001.318} This argument shows that for any
torsion class $\TC\in\Tsn(H^3(X,\bbZ))$ there is a \emph{smooth} Azumaya
bundle with Dixmier--Doaudy invariant $\TC.$ \end{rem*}

\subsection{Twisted $K$--theory}

Let $X$ be a manifold and let $\cJ$ be a locally trivial bundle of algebras
over $X$ with fibre $\cK$ and structure group $\Aut(\cK).$ Two such bundles
are isomorphic if and only if they have the same Dixmier--Douady invariant
$\delta(\cJ)\in H^3(X,\bbZ).$ The \emph{twisted $K$--theory} of $X$ (with
compact supports) has been
defined by Rosenberg \cite{Rosenberg1} as
\begin{equation}
\begin{gathered}
\Kc^j(X,\cJ)=K_j(C_0(X,\cJ))\quad j=0,1,
\end{gathered}
\label{twisted k-theory}\end{equation}
where and $K_\bullet(C_0(X,\cJ))$ denotes the topological $K$--theory of the
\Cstar--algebra of continuous sections of $\cJ$ that vanish outside a
compact subset of $X.$ In case $X$ is compact we use the notation
$K^j(X,\cJ).$ The space $K^j(X,\cJ)$ or $\Kc^j(X,\cJ)$ is an abelian
group. It is tempting to think of the twisted $K$--theory of $X$ as determined
by the class $\TC=\delta (\cJ).$ However, this is not strictly speaking
correct; whilst it is the case that any other choice $\cJ'$ such that
$\delta(\cJ') = \TC$ is isomorphic to $\cJ$ and therefore there is an
isomorphism $K^j(X,\cJ)\cong K^j(X;\cJ')$ this isomorphism itself is not
unique, nor is its homotopy class. However both the abelian group structure
and the module structure on $K^0(X;\cJ)$ over $K^0(X)$ arising from
tensor product are natural, ie, are preserved by such isomorphisms. With
this caveat one can use the notation $K^j(X,\TC)$ to denote the twisted
$K$--theory with Dixmier--Douady invariant $\TC\in H^3(X,\bbZ).$

In the case of principal interest here, when $\TC\in\Tsn(H^3(X,\bbZ))$ we
will take $\cJ=\cK_{\cA}=\cA\otimes\cK$ where $\cA$ is an Azumaya bundle
and use the notation
\begin{equation}
K^j_c(X,\cA)=K_j(C_0(X,\cA\otimes\cK))=K_j(C_0(X,\cJ))\quad j=0,1
\label{mms2001.357}\end{equation}
and $K^j(X,\cA)$ in the compact case. As noted above, if $\cA$ and $\cA'$
are two Azumaya bundles with the same Dixmier--Douady invariant, then as
bundles of algebras, $\cA\otimes\End(V)\equiv\cA'\otimes\End(W)$ for some
vector bundles $V$ and $W.$ Whilst there is no natural isomorphism of
$\End(V)\otimes\cJ$ and $\cJ,$ these bundles of algebras are
isomorphic. Moreover $\Pi_0$ of the group induced by diffeomorphisms on
$\cJ$ is naturally isomorphic to $H^2(X,\bbZ).$ It follows that the
isomorphism between $K^0(X,\cA)$ and $K^0(X,\cA')$ is determined up to the
action of the image of $H^2(X,\bbZ),$ as isomorphism classes, in $K^0(X)$
acting on $K^0(X,\cA)$ through the module structure.

There are alternate descriptions of $K^0(X,\cA).$ A description in terms of
the twisted index map is mentioned in \cite{Rosenberg1} and we give a
complete proof here. Let $Y_{\cA}$ be the principal $\PU(\H)= \Aut(\cK)$
bundle over $X$ associated to $\K_\cA=\cA\otimes\K$ and let $\Fred_\cA =
\left(Y_\cA \times \Fred(\cH) \right)/\PU(\cH)$ be the bundle of twisted
Fredholm operators where $\Fred(\cH)$ denotes the space of Fredholm 
operators on
$\cH.$

We shall use the short exact sequence of \Cstar--algebras,
$$
0\longrightarrow\K\longrightarrow \B \longrightarrow \B/\K\longrightarrow
0
$$
where $\B/\K$ is the Calkin algebra. It gives rise to the short exact
sequence of \Cstar--algebras of sections,
\begin{equation}
0\longrightarrow C(X, \K_{\cA})\longrightarrow C(X, \B_\cA)
\longrightarrow C(X, (\B/\K)_\cA)\longrightarrow0
\label{exact1}\end{equation}
where the bundles $\B_\cA$ and $(\B/\K)_\cA$ are also associated
to $Y_\cA$.

\begin{lemma}\label{Bundle_realization} Let $X$ be a compact manifold and
$\cA$ an Azumaya bundle over $X$. Then the {twisted $K$--theory}
$K^0(X,\cA)$ is isomorphic to the Grothendieck group of Murray--von Neumann
equivalence classes of projections in $C(X, \K_\cA).$ Similarly,
$K^0(X,\B_\cA)$ is isomorphic to the Grothendieck group of Murray--von
Neumann equivalence classes of finite rank projections in $C(X, \B_\cA).$
\end{lemma}

\begin{proof} This is a consequence of Lemma~\ref{stably_unital} above,
which asserts that $C(X,K_\cA)$ is stably unital, and of Proposition
5.5.5. in \cite{Dixmier-Douady1}.
\end{proof}


\begin{proposition}\label{Index-bundle}
For any compact maniold $X$, there is a twisted index map, defined
explicitly in \eqref{indexp} below, giving an isomorphism.
\begin{equation}
\index \co  \pi_0(C(X, \Fred_\cA)) \stackrel{\sim}{\longrightarrow} K^0(X,
\cA).
\label{Twisted_index}\end{equation}
\end{proposition}

\begin{proof} First we show that $K_0 (C(X,\B_\cA)) =\{0\}.$ We remind the
reader that
\begin{equation*}
\GL(n,\B)\simeq \GL(\B)
\end{equation*}
is connected, so $K_1(\B)=0.$ Also $K_0(\B)=0$ because all infinite
projections in $\B$ are equivalent to $I$, so all finite projections are
equivalent to $0.$

A finite projection in $C(X,\B_\cA)$ is a cross section of finite
projections, $x \longmapsto P_{x}, x \in X;$ see
Lemma~\ref{Bundle_realization}. The set of partial isometries at $x$
connecting $I$ to $I-P_x$ is a principal homogeneous space for the full
unitary group and hence contractible.  A cross section of partial
isometries makes $I$ equivalent to $I-P$ and $P$ equivalent to $0.$ Hence
$K_0 (C(X,\B_\cA)) = {0}.$

We leave to the reader the proof that $K_1 (C(X,\B_\cA))=\{0\},$ using the
fact that $GL(\B)$ is contractible.

Now consider the six term exact sequence in
$K$--theory
\begin{equation}
\begin{CD}
K^1(X, \cA) @>>> K_1(C(X, \B_\cA)) @>>>
K_1(C(X, (\B/\K)_\cA)) \\
                @AAA        &               &
@VV\index' V \\
K_0(C(X, (\B/\K)_\cA))   @<<< K_0(C(X, \B_\cA)) @<<<
K^0(X, \cA)
\end{CD}
\label{kexact}\end{equation}
arising from \eqref{Twisted_index}.

From this, using \eqref{kexact}, we get
$K_1 (C(X,\B_X)\cA) \simeq K^0 \left( {X,\cA} \right).$ However, note that
by definition
$$
\begin{array}{lcl}
K_1(C(X, (\B/\K)_\cA)) &=&\displaystyle\lim_{\to} \;\pi_0\left( \GL(n,
C(X, (\B/\K)_\cA))\right)\\[+7pt]
&=& \displaystyle\lim_{\to} \;\pi_0(C(X, \GL(n,\B/\K)_\cA))
\end{array}
$$
where $\GL(n, A)$
denotes the group of invertible $n\times n$ matrices with entries
in the \Cstar--algebra $A$. In the case of the Calkin algebra,
$\B(\H) \otimes M_n(\C) \cong \B(\H\otimes \C^n)$ and
$\K(\H) \otimes M_n(\C) \cong \K(\H\otimes \C^n)$ from which it follows
that
$$
\B(\H)/\K(\H)\otimes M_n(\C) \cong \B(\H\otimes \C^n)/\K(\H\otimes \C^n).
$$
Therefore $\GL(n, \B(\H)/\K(\H)) \cong \GL(1, \B(\H\otimes
\C^n)/\K(\H\otimes \C^n)),$
so
$$
K_1(C(X, (\B/\K)_\cA)) \cong \pi_0(C(X, \GL(1,\B/\K)_\cA)),
$$
and we obtain the isomorphism
\begin{equation}
\index' \co  \pi_0(C(X, \GL(1, \B/\K)_\cA))\longrightarrow K^0(X,\cA).
\label{indexp}\end{equation}
\end{proof}

Thus when $\TC$ is a torsion class, the corresponding twisted $K$--theory can
be described just as if the \Cstar--algebra $C(X, \K_\cA)$ had a unit. In
this case the index isomorphism can also be written in familiar form.

\begin{proposition}\label{Index-bundle-old} Suppose $X$ is a compact manifold
and $\cA$ is an Azumaya bundle over $X.$ Then given a section $s\in C(X,
\Fred_\cA)$ there is a section $t\in C(X, \K_\cA)$ such that $\index(s+t)
= p_1 - p_2$, where $p_1, p_2$ are projections in $C(X, \K_\cA)$
representing the projection onto the kernel of $(s+t)$ and the projection
onto the kernel of the adjoint $(s+t)^*$ respectively. Moreover the map
$$
\begin{array}{rcl}
\pi_0(C(X, \Fred_\cA))&\longrightarrow & K^0(X, \cA)\\[+7pt]
[[s]] &\longrightarrow  & [\index(s+t)]
\end{array}
$$
is well defined and is independent of the choice of $t\in C(X,\K_\cA).$
\end{proposition}

\subsection{Projective vector bundles}

Consider the short exact sequence of groups
\begin{equation}
\bbZ_n\longrightarrow \SU(n)\overset\pi\longrightarrow \PU(n),\ n\in\bbN
\label{mms2001.358}\end{equation}
and the associated (determinant) line bundle $L$ over $\PU(n).$ The
fiber at $p\in\PU(n)$ is
\begin{multline}
L_p=\{(a,z)\in\SU(n)\times\bbC;\pi(a)=p\}/\sim,\\
(a,z)\sim(a',z')\text{ if }a'=ta,\ z'=tz,\ t\in\bbZ_n.
\label{mms2001.359}\end{multline}
This is a primitive line bundle over $\PU(n)$ in the sense that there is a
natural $\SU(n)$ action on the total space of $L,$
\begin{equation}
l_a\co (g,z)\longmapsto (ga,z)
\label{mms2001.360}\end{equation}
which induces a natural isomorphism
\begin{equation}
L_{pq}\equiv L_p\otimes L_q\ \forall\ p,q\in\PU(n).
\label{mms2001.361}\end{equation}
Indeed, if $a\in\SU(n)$ and $\pi(a)=p$ then $l_q\co L_{b}\longrightarrow L_{ab}$
and $l_a\co L_{\Id}\equiv\bbC\longrightarrow L_q$ combine to give an isomorphism
\eqref{mms2001.361} which is independent of choices. From the definition
there is an injection
\begin{equation}
i_q\co \pi^{-1}(q)\hookrightarrow L_q,\ q\in\PU(n)
\label{mms2001.362}\end{equation}
mapping $a$ to the equivalence class of $(a,1).$ Thus $a\in\SU(n)$ fixes a
trivialization $e_a\co L_{\pi(a)}\longrightarrow \bbC,$ determined by
$(a,1)\longmapsto 1.$

Let $P=P({\mathcal{A}})$ be the $\PU(n)$ bundle associated to an Azumaya
bundle. Thus the fiber $P_x$ at $x\in X$ consists of the algebra
isomomorphisms of $\cA_x$ to $M(n, \bbC).$ By a \emph{projective} vector
bundle over $X,$ associated to $\cA,$ we shall mean a complex vector bundle
$E$ over $P({\cA})$ with a smooth family of linear isomorphisms
\begin{equation}
\gamma_p\co p^*E\longrightarrow E\otimes\fL_p,\ p\in\PU(n)
\label{mms2001.349}\end{equation}
satisfying the compatibility condition
\begin{equation}
\gamma _{pp'}=\gamma_{p}\circ \gamma _{p'}
\label{mms2001.350}\end{equation}
in the sense that on the right $\gamma_p\co p^*E\otimes\fL_{p'}\longrightarrow
E\otimes\fL_{p}\otimes\fL_{p'}\longrightarrow E\otimes\fL_{pp'}$ using
\eqref{mms2001.361}. In fact $\gamma$ lifts to an action of $\SU(n)$:
\begin{equation}
\tilde\gamma_a\co p^*E\longrightarrow E,\ \tilde\gamma_a=(\Id\otimes
e_a)\gamma_{\pi(a)}.
\label{mms2001.363}\end{equation}
Thus projective vector bundles are just a special case of $\SU(n)$--equivariant
vector bundles over $P.$ This has also been studied in the case when
$P$ is a bundle gerbe \cite{Murray1} in which case $E$ is known as a
bundle gerbe module \cite{Bouwknegt-Carey-Mathai-Murray-Stevenson1}.

A bundle homomorphism between two projective bundles $E$ and $F$ is itself
projective if it intertwines the corresponding isomorphisms
\eqref{mms2001.349}. Since there is a natural isomorphism
$\hom(E_m,F_m)\equiv\hom(E_m\otimes L,F_m\otimes L)$ for any complex line
$L,$ the identifications $\gamma_p$ act by conjugation on $\hom(E,F)$ and
give it the structure of a $\PU(n)$--invariant bundle. Thus the invariant
sections, the projective homomorphisms, are the sections of a bundle
$\hom_{\proj}(E,F)$ over $X.$

Just as a family of finite rank projections, forming a section of
$\mathcal{C}(X,\mathcal{K}),$ fixes a vector bundle over $X,$ so a section
of a model twisted bundle $\mathcal{A}\otimes\mathcal{K},$ with values
in the projections, fixes a projective vector bundle over $X.$ Thus, if
$m\in P_x$ then by definition $m\co \cA_x\longrightarrow {M}(n, \bbC)$ is an
algebra isomorphism. The projection $\mu_x\in\cA_x\otimes~\!\cK$ thus becomes
a finite rank projection in ${M}(n, \bbC)\otimes\cK.$ Using a fixed
identification
\begin{equation}
R\co \bbC^n\otimes\cH\longrightarrow \cH
\label{mms2001.364}\end{equation}
and the induced identification $\Ad(R)$ of ${M}_n(\bbC)\otimes\cK$ and
$\cK(\cH),$ this projection may be identified with its range
$E_m\subset\cH.$ The continuity of this operation shows that the $E_m$ form
a vector bundle $E$ over $P.$ To see that $E$ is a projective
vector bundle observe that under the action of $p\in\PU(n)$ on $P_x,$
replacing $m$ by $mp,$ $E_m$ is transformed to $E_{m'}=R(a\otimes\Id)E_m$
where $a\in\SU(n)$ is a lift of $p,$ $\pi(a)=p.$ Using the choice of $a$ to
trivialize $L_p,$ the resulting linear map
\begin{equation}
\gamma_p\co E_m\longrightarrow E_{m'}\otimes L_p^{-1}
\label{mms2001.365}\end{equation}
is independent of choices and satisfies \eqref{mms2001.350}.

The direct sum of two projective vector bundles over $P$ is again a
projective bundle so there is an associated Grothendieck group of
projective $K$--theory over $P,$ $K^0_{\proj}(P).$ The discussion above of the
equivalence of projections in $C(X;\cA\otimes\cK)$ and projective vector
bundles then shows the natural equality of the corresponding Grothendieck
groups, just as in the untwisted case.

\begin{lemma}\label{mms2001.352} If $P$ is the principal $\PU(n)$ bundle
over $X$ associated to an Azumaya bundle, then $K^0_{\proj}(P)$ is
canonically isomorphic to $K^0(X;\mathcal{A}).$
\end{lemma}

Similar conclusions hold for $K$--theory with compact supports if the base is
not compact; we denote the twisted $K$--groups with compact support
$\Kc(X,\mathcal{A}).$

A projective vector bundle may be specified by local trivializations
relative to a trivialization of $\mathcal{A}$ and we proceed to discuss the
smoothness and equivalence of such trivializations.

Consider a `full' local trivialization of the Azumaya bundle
$\mathcal{A}$ over a good open cover $\{U_a\}_{a\in A}$ of the base $X.$ Thus,
there are algebra isomorphisms
\begin{equation*}
F_a\co \mathcal{A}\big|_{U_a}\longrightarrow U_a\times M(n,\bbC)
\label{mms2001.346}\end{equation*}
with lifted transition maps, chosen to be to continuous (or smooth if the
$G_{ab}$ are smooth)
\begin{equation}
G_{ab}\co U_{ab}=U_a\cap U_b\longrightarrow \SU(n)\text{ such that }
G_{ab}\equiv F_a\circ F_b^{-1}\Min\PU(n).
\label{mms2001.347}\end{equation}
That is, the transition maps for $\mathcal{A}$ over $U_{ab}$ are given by
the adjoint action of the $G_{ab}.$ Thus the Dixmier--Douady cocycle
associated to the trivialization is
\begin{equation}
\TZ_{abc}=G_{ab}G_{bc}G_{ca}\co U_{abc}\longrightarrow
\bbZ_n\subset {\rm U}(1)\subset\bbC^*.
\label{mms2001.348}\end{equation}
Such a choice of full local trivialization necessarily gives a local
trivialization of the associated $\PU(n)$ bundle, $P,$ and also gives a local
trivialization of the determinant bundle $L$ over $P.$
The Dixmier--Douady class of $\cA$ is defined as the cohomology class 
$[\theta] \in H^2(X, \underline{{\rm U}(1)})
\cong H^3(X, \bbZ)$. Whereas it is true that one can pick
a ``constant'' $\bbZ_n$--valued, 2--cocycle that represents the 
Dixmier--Douady class as shown above, the equivalence
relation (ie cohomology class) is not in $H^2(X, \bbZ_n)$, but rather in the sheaf
cohomology  $H^2(X, \underline{{\rm U}(1)})$ (which is isomorphic to $H^3(X, \bbZ)$).

\emph{Projective vector bundle data,} associated to such a full
trivialization of $\cA,$ consists of complex vector bundles $E_a,$ of some
fixed rank $k,$ and transition maps $Q_{ab}\co E_b\longrightarrow E_a$ over
each $U_{ab}$ satisfying the weak cocycle condition
\begin{equation}
Q_{ab}Q_{bc}=\TZ_{abc}Q_{ac}
\label{mms2001.336}\end{equation}
where $\TZ$ is given by \eqref{mms2001.348}. Two sets of such data $E_a,$
$Q_{ab}$ and $E'_a,$ $Q'_{ab}$ over the same cover are \emph{equivalent} if
there are bundle isomorphisms $T_a\co E_a\longrightarrow E_a'$ such that
$Q'_{ab}=T_aQ_{ab}T_b^{-1}$ over each $U_{ab}.$

Associated with the trivialization of $\mathcal{A}$ there is a particular
set of projective vector bundle data given by the trivial bundles $\bbC^n$
over the $U_a$ and the transition maps $G_{ab}.$ We will denote this data
as $E_\tau$ where $\tau$ denotes the trivialization of $\mathcal{A}.$ The
Azumaya algebra $\cA$ may then be identified with
$\hom_{\proj}(E_\tau,E_\tau)$ for any of the projective vector bundles
$E_\tau.$

\begin{lemma}\label{mms2001.321} Projective vector bundle data with respect
to a full trivialization of an Azumaya bundle lifts to define a projective
vector bundle over the associated principal $\PU(n)$ bundle; all projective
bundles arise this way and projective isomorphisms of projective bundles
corresponds to equivalence of the projective vector bundle data.
\end{lemma}

\begin{proof} The given trivialization of $\mathcal{A},$ over each $U_a$
defines a section of $s_a\co U_a$ $\longrightarrow P=P({\mathcal{A}})$ over
$U_a.$ Using this section we may lift the bundle $E_a$ to the image of the
section and then extend it to a bundle $E^{(a)}$ on the whole of
$P\big|_{U_a}$ which is projective, namely by setting
\begin{equation}
E{(a)}_{s_a(x)p}=E_a(x)\otimes L_p\ \forall\ p\in\PU(n),\ x\in U_a
\label{mms2001.351}\end{equation}
and taking the ismorphism $\gamma_p$ over $U_a$ to be given by the identity
on $E_a(x).$ Over each intersection we then have an isomorphism
$$
Q_{ab}(x)\co E_a(x)=E^{(a)}(s_a(x))\longrightarrow E^{(b)}(s_b(x))=E_b(x).
$$
Now, from the trivialization of $P$ we have $s_b(x)=s_a(x)g_{ab}(x)$ where
$G_{ab}\co U_{ab}$ $\longrightarrow \SU(n)$ and $g_{ab}$ is the projection of
$G_{ab}$ into $\PU(n).$ The choice of $G_{ab}$ therefore also fixes a
trivialization of the determinant bundle
$L_{g_{ab}(x)}\longrightarrow\bbC.$ Since $E^{(b)}(s_b(x))$ is identified
with $E^{(b)}(s_a(x))\otimes L_{g_{ab}(x)}^{-1}$ by the primitivity, this
allows $Q_{ab}$ to be interpreted as the transition map from $E^{(a)}$ to
$E^{(b)}$ over the preimage of $U_{ab}.$ Furthermore the weak cocycle
condition \eqref{mms2001.336} now becomes the cocycle condition
guaranteeing that the $E^{(a)}$ combine to a globally defined, projective,
bundle over $P.$

This argument can be reversed to construct projective vector bundle data
from a projective vector bundle over $P$ and a similar argument shows that
projective bundle isomorphisms correspond to isomorphisms of the projective
vector bundle data.
\end{proof}

Thus we may simply describe `projective vector bundle data' as a local
trivialization of the corresponding projective bundle, where this also
involves the choice of a full local trivialization of $\cA.$ The projective
vector bundle data $E_\tau$ associated to a full trivialization of $\cA$
thus determines a projective vector bundle, which we may also denote by
$E_\tau,$ over $P(\cA).$ In particular there are projective vector bundles of
arbitrarily large rank over $P.$

Many of the standard results relating $K$--theory to vector bundles carry over
to the twisted case. We recall two of these which are important for the
proof of the index theorem.  Lemma \ref{mms2001.327} below essentially 
follows from definitions, the morphism
$C_0(U, \cA) \to C(M, \cA)$, and the functoriality of
$K$--theory.

\begin{lemma}\label{mms2001.327} If $U\subset X$ is an open set of a
compact manifold then for any Azumaya bundle $\mathcal{A},$ over $X,$ there
is an extension map
\begin{equation}
K^0_c(U,\mathcal{A}_U)\longrightarrow K^0(X,\mathcal{A}).
\label{mms2001.328}\end{equation}
\end{lemma}

\begin{proof} This proceeds in essentially the usual way. An element of
$K^0_c(U,\mathcal{A}_U)$ is represented in terms of a full local
trivialization of $\mathcal{A}$ by a pair of sets of projective vector data
$E_1,$ $E_2$ over $U$ and a given bundle isomorphism between them outside a
compact set, $c\co E_1\longrightarrow E_2$ over $U\setminus K.$ If $E_\tau$ is
the projective vector bundle data associated to the trivialization of
$\mathcal{A}$ then we may embed $E_2$ as a projective subbundle of
$E_\tau^q$ for some integer $q.$ To see this, first choose an embedding $e$
of $E_2$ as a subbundle of $E_\tau^l$ over $P$ for some large $l.$ Then
choose a full trivialization of $\cA$ and a partition of unity $\psi _a$
subordinate on $X$ to the open $U_a.$ Over the premimage of each $U_a,$
$e_a$ can be extended uniquely to a projective embedding $e_a$ of $E$ in
$E_\tau^l.$ The global map formed by the sum over $a$ of the $\psi _ae_a$
gives a projective embedding into $E_\tau^{q}$ where $q=lN$ and $N$ is the
number of sets in the open cover. Then $E_\tau^q=E_2\oplus F$ for some
complementary projective bundle $F$ over $U.$ The data $E_1\oplus F$
over $U$ and
$E_\tau^q$ over $Z\setminus K$ with the isomorphism $c\oplus\Id_F$ over
$U\setminus K$ determine a projective vector bundle over $Z.$ The element
of $K^0(X;\mathcal{A})$ represented by the pair consisting of this bundle
and $E_\tau^q$ is independent of choices, so defines the extension map.
\end{proof}

\begin{proposition}\label{mms2001.319} For any real vector bundle
$\pi\co V\longrightarrow X$ and Azumaya bundle $\mathcal{A}$ over $X$ a
section over the sphere bundle, $SV,$ of $V$ of the isomorphism bundle of
the lifts of two projective vector bundles over $X,$ determines an element of
$K^0_c(V,\pi^*\mathcal{A});$ all elements arise this way and two
isomorphisms give the same element if they are homotopic after
stabilization with the identity isomorphisms of projective bundles.
\end{proposition}

\begin{proof} The proof is the same as in the untwisted case.
\end{proof}

If $W$ is a vector bundle over $X$ and $E$ is a projective bundle over a
$\PU(n)$ bundle over $X$ then $E\otimes\pi^*W$ is naturally a projective
bundle. This operation extends to make $K^0_{\proj}(P)$ a module over
(untwisted) $K^0(X)$ and hence, in view of Lemma~\ref{mms2001.352} gives a
module structure
\begin{equation}
K^0(X)\times K^0(X;\mathcal{A})\longrightarrow K^0(X;\mathcal{A}).
\label{mms2001.353}\end{equation}

\subsection{The Chern character of projective vector bundles}\label{SS.Chern}

The Chern character
\begin{equation}
\Ch_\A\co K^0(X;\mathcal{A})\longrightarrow H^{\ev}(X;\bbQ)
\label{mms2001.354}\end{equation}
may be defined in one of several equivalent ways.  
It is known that $K^0(X)$ and $K^0(X, \cA)$ are isomorphic
after tensoring with $\bbQ$,  \cite{MR39:5586a} (see also \cite{Witten3}, \cite{Nistor-Troisky1}). 
Since the Chern Character on $K^0(X)$ factors through $K^0(X)\otimes\bbC$ it
is also defined in the twisted case.

To define \eqref{mms2001.354} directly using the Chern--Weil approach we
note that the local constancy of the C\v ech 2--cocycle $\TZ$ in
\eqref{mms2001.336} allows a connection to be defined directly on such
projective vector bundle data. That is, despite the failure of the usual
cocycle condition, there exist connections $\nabla^a$ on each of the
bundles $E_a$ which are identified by the $Q_{ab}.$ To see this, simply
take arbitary connections $\tilde\nabla^a$ on each of the $E_a$ and a
partition of unity $\phi _a$ subordinate to the cover. Now define a new
connection on $E_a$ by
\begin{equation*}
\nabla^a=\phi _a\nabla^a+\sum\limits_{b\not=a}\phi _bQ_{ab}^*\nabla^b.
\label{mms2001.355}\end{equation*}
These are consistent under the transition maps. Thus the curvature of this
collective connection is a well-defined section of the endomorphism bundle
of the given projective vector bundle data. As such the usual Chern--Weil
arguments apply and give the Chern character \eqref{mms2001.354}. Lifted to
$P$ this connection gives a projective connection on the lift of $E$ to a
projective bundle; this also allows the Chern character to be defined
directly on $K^0_{\proj}(P).$

Either of these approaches to the Chern character show that it distributes
over the usual Chern character on $K^0(X)$ under the action
\eqref{mms2001.353}. In particular it follows that \eqref{mms2001.354} is
an isomorphism over $\bbQ.$

We may also define the Chern character by reducing to the standard case by
taking tensor powers. Thus if $E$ is a projective vector bundle associated
to the Azumaya bundle $\cA$ and $n$ is the order of the
associated Dixmier--Douady invariant then
$$
\phi^*E^{\otimes n} \cong \pi_P^*E^{\otimes n}
$$
and therefore $E^{\otimes n} = \pi^*(F)$ for some vector bundle
$F\to X$.

Observe that the Chern character of $E$ satisfies,
\begin{equation} \label{ind1}
\Ch(E^{\otimes n}) = \Ch(E)^n = \pi^*(\Ch(F)).
\end{equation}
We claim that $\Ch(E) = \pi^*(\Lambda)$ for some cohomology class $\Lambda\in
H^{even}(X, \Q).$ First observe that \eqref{ind1} implies that the degree
zero term, which is a constant term, is of the desired form. Next assume
that the degree $2k$ component of the Chern character satisfies $\Ch_k(E) =
\pi^*(\Lambda_k)$ for some
cohomology class $\Lambda_k\in H^{2k}(X, \Q)$. Then \eqref{ind1} implies that
the degree (2k+2) component
$\Ch_{k+1}(E^{\otimes n}) = \pi^*(\Ch_{k+1}(F))$.
But the left hand side is of the form of the Chern character,
\begin{equation} \label{chernind1}
\Ch_{k+1}(E^{\otimes n}) = a_0 \Ch_{k+1}(E) +
\sum_{|I|=k+1, r>1} a_I \Ch_{I}(E)
\end{equation}
where $I = (i_1, \ldots, i_r)$, $|I|= i_1+\cdots + i_r$,
$\Ch_{I}(E) = \Ch_{i_1}(E)\cup\cdots \cup \Ch_{i_r}(E)$
and $a_0, a_I \in \mathbb Q$ are such that $a_0\ne 0$. By the
induction hypothesis, we deduce that $\Ch_{k+1}(E)$ is of the
form $\pi^*(\Lambda_{k+1})$ for some
cohomology class $\Lambda_{k+1}\in H^{2k+2}(X, \Q)$.
This proves the claim.

Then the Chern character of the projective vector bundle $E$ above is given by
\begin{equation} \label{cherncharacter}
\Ch_\A(E) = \Lambda \in H^{even}(X, \Q).
\end{equation}
That is, the lift of the Chern character of the projective vector bundle
$E$ to $P$ coincides with the ordinary Chern character of $E.$  The
following properties of the Chern character of projective vector bundles
follow from the corresponding properties of the Chern character of vector
bundles.

\begin{lemma}\label{chern} Let $\A$ be an Azumaya bundle over $X$ and
$P$ be the principal ${\rm PU}(n)$ bundle associated to $\A.$ Let $E\to P$
be a projective vector bundle over $X$, associated to $\A$, then the Chern
character defined as in \eqref{cherncharacter} above has the
following properties.
\begin{enumerate}
\item
If $E'\to X$ is another projective vector bundle, then
$$
\Ch_\A(E\oplus E') = \Ch_\A(E) + \Ch_\A(E'),
$$
so the Chern character is a homomorphism
$$
\Ch_\A\co  K^0(X, \A) \to H^{even}(X, \Q).
$$
\item
The degree 2 component of the Chern character $\Ch_\A(E)$ coincides with the
first Chern class of the determinant line bundle $\det(E)\to X$.
\end{enumerate}
\end{lemma}

\section{The analytic index}\label{Analytic}

For two ordinary vector bundles, $E^\pm,$ over a compact manifold $Z,$ the
space $\Psi^m(Z;E^+,E^-)$ of pseudodifferential operators of order $m$
mapping $\Cinf(Z;E^+)$ to $\Cinf(Z;E^-)$ may be identified
naturally with the tensor product
\begin{equation}
\Psi^m(Z;E^+,E^-)=\Psi^m(Z)\otimes_{\Cinf(Z^2)}\Cinf(Z^2;\Hom(E^+,E^-))
\label{mms2001.307}\end{equation}
where $\Hom(E^+,E^-)$ is the `big' homomorphism bundle over $Z^2$ which has
fiber $\hom(E^+_z,E^-_{z'})$ at $(z',z)$ and $\Psi^m(Z)$ is the space of
pseudodifferential operators acting on functions. The latter is a module
over $\Cinf(Z^2)$ through its realization as a space of Schwartz
kernels. In particular
\begin{equation*}
\Psi^m(Z;E,E)=\Psi^m(Z;E)=\Psi^m(Z)\otimes_{\Cinf(Z^2)}\Cinf(Z^2;\Hom(E))
\label{mms2001.308}\end{equation*}
when the two bundles coincide.

For a fibration $\phi\co M\longrightarrow X,$ with compact boundaryless
fibres, the bundle of pseudodifferential operators acting on sections of
vector bundles over the total space may be similarly defined. Note that the
operators act fibre-wise and so commute with multiplication by functions on
the base. If $\Hom(M^2_\phi,E^+,E^-)$ denotes the bundle over the fibre
product, which is the `big' homomorphism bundle on each fibre, then again
\begin{equation*}
\Psi^m(M/X;E^+,E^-)=\Psi^m(M/X)\otimes_{\Cinf(M^2_\phi)}
\Cinf(M^2_\phi;\Hom(E^+,E^-))
\label{mms2001.311}\end{equation*}
is the bundle of operators to which the usual families index theorem
applies. Here $\Psi^m(M/X)$ is the bundle of pseudodifferential operators
acting on functions on the fibres.

Now, let $\mathcal{A}$ be an Azumaya bundle over the base of the fibration
$\phi.$ Consider a projective vector bundle, $E,$ over the lift to $M$
of the principal $\PU(N)$ bundle associated to $\mathcal{A}.$ Given
a local trivialization of $\mathcal{A}$ there is a bundle trivialization of
$E$ with respect to the lift of the trivialization to $M.$ We shall call
this a \emph{basic bundle trivialization.}

\begin{lemma}\label{mms2001.312} If $E^+$ and $E^-$ are two projective
       vector bundles over the lift to $M$ of the $\PU(N)$ bundle associated
       to an Azumaya bundle on $X$ then the big homomorphism
       bundles $\Hom(Q^+_a,Q^-_a),$ arising from basic bundle trivializations
       of the $E^\pm$ define a vector bundle $\Hom(E^+,E^-)$ over $M^2_\phi.$
\end{lemma}

\begin{proof} As already noted, a local trivialization of $\mathcal{A}$
       over the base gives a trivialization of the associated $\PU(N)$ bundle,
       and hence of its lift to $M.$ This leads to bundle trivializations
       $Q^\pm_a,$ over the elements $\phi ^{-1}(U_a)$ of this open cover, of
       the projective bundles $E^\pm.$ Since the transition maps act by the
       adjoint action, the scalar factors cancel and the `big' homomorphism
       bundles between the $Q_a^\pm$ now patch to give a global bundle
       $\Hom(E^+,E^-)$ over $M^2_\phi.$
\end{proof}

If $E_i,$ $i=1,2,3,$ are three such projective bundles then, just as for
the usual homomorphism bundles, there is a bilinear product map
\begin{equation}
\Hom(E_1,E_2)\otimes\Hom(E_2,E_3)\big|_C\longrightarrow \psi^*\Hom(E_1,E_2)
\label{mms2001.316}\end{equation}
where $C$ is the central fiber diagonal in $M^2_\phi\times M^2_\phi$ and
$\psi\co C\longrightarrow M^2_\phi$  is projection off the middle factor. This
reduces to the composition law for
$\hom(E_i,E_j)=\Hom(E_i,E_j)\big|_{\Diag}$ on the diagonal.

We may now simply define the algebra of twisted (fiber-wise) pseudodifferential
operators as
\begin{equation}
\Psi^m(M/X;E^+,E^-)=\Psi^m(M/X)\otimes_{\Cinf(M^2_\phi)}
\Cinf(M^2_\phi;\Hom(E^+,E^-)).
\label{mms2001.309}\end{equation}
Restricted to open sets in the base over which $\mathcal{A}$ is
trivialized, this reduces to the standard definition. Thus, the symbol
sequence remains exact
\begin{multline}
0\longrightarrow \Psi^{m-1}(M/X;E^+,E^-)\hookrightarrow\\
\Psi^m(M/X;E^+,E^-)\overset{\sigma _m}{\longrightarrow}
S^{[m]}(S^*(M/X);\phi^*\hom(E^+,E^-))\longrightarrow 0
\label{mms2001.310}\end{multline}
with the proof essentially unchanged. Here
$S^{[m]}(S^*(M/X);\rho^*\hom(E^+,E^-))$ is the quotient space of symbolic
sections of order $m,$ by symbolic sections of order $m-1,$ of $\rho
^*\hom(E^+,E^-)$ as a bundle over $S^*(M/X),$ the fibre cosphere
bundle, $\rho\co S(M/X)\longrightarrow M$ being the projection. Similarly the
usual composition properties carry over to this
twisted case, since they apply to the local families. For any three
projective vector bundles $E_i,$ $i=1,2,3,$ over the lift of the same
$\PU(N)$ bundle from the base
\begin{multline}
\Psi^m(M/X;E_2,E_3)\circ\Psi^{m'}(M/X;E_1,E_2)\subset\Psi^{m+m'}(M/X;E_1,E_3),\\
\sigma _{m+m'}(A\circ B)=\sigma _m(A)\circ\sigma _{m'}(B).
\label{mms2001.315}\end{multline}.

For any basic bundle trivialization of a projective vector bundle with
respect to a local trivialization of $\mathcal{A}$ the spaces of sections
of the local bundles form infinite-dimensional projective bundle data over
the base, associated to the same trivialization of $\mathcal{A}.$ More
generally, for any fixed real number, $m,$ the spaces of
Sobolev sections of order $m$ over the fibres form projective Hilbert
bundle data over the base; we will denote the corresponding projective bundle
$H^m(M/X;E).$ The boundedness of pseudodifferential operators on Sobolev
spaces then shows that any $A\in\Psi^m(M/X;E^+,E^-)$ defines a bounded
operator
\begin{equation}
A\co H^{m_1}(M/X,E^+)\longrightarrow H^{m_2}(M/X;E^-)\text{ provided
}m_1\ge m_2+m.
\label{mms2001.317}\end{equation}
If $m_1>m_2+m$ this operator is compact.

It is possible to choose quantization maps as in the untwisted case. To do
so, choose basic bundle trivializations and quantization maps, that is
right inverses for the symbol map, for the local bundles $Q^\pm_a.$ Using a
partition of unity on the base this gives a global quantization map:
\begin{multline}
q_m\co S^{[m]}(S^*(M/X);\rho^*\hom(E^+,E^-))\longrightarrow \Psi^m(M/X;E^+,E^-),\\
\sigma _m\circ q_m=\Id,\
q_m\circ\sigma_m-\Id\co \Psi^m(M/X;E^+,E^-)\longrightarrow
\Psi^{m-1}(M/X;E^+,E^-).
\label{mms2001.314}\end{multline}
By definition a projective family in $\Psi^m(M/X;E^+,E^-)$ is elliptic if
$\sigma _m$ is invertible, with inverse in
$S^{[m]}(S^*(M/X);\rho^*\hom(E^-,E^+)).$ Directly from the symbolic
properties of the algebra, this is equivalent to there being a parameterix
$B\in\Psi^{-m}(M/X;E^-,E^+)$ such that $A\circ B-\Id\in\Psi^{-1}(M/X;E^-)$
and $B\circ A-\Id\in\Psi^{-1}(M/X;E^+).$ These `error terms' give compact
maps, for $m_1=m_2+m,$ in \eqref{mms2001.317}. Thus the elliptic family
consists of Fredholm operators. It follows from the discussion in
Section~\ref{Review} that the family defines a twisted $K$--class using
\eqref{Twisted_index}. To see this class more concretely, as in the
untwisted case, we may perturb the family so that the index bundle gives
projective vector bundle data with respect to the given trivialization of
$\mathcal{A}.$ Locally in the base a bundle map from an auxilliary vector
bundle, over the base, may be added to make the family surjective. Choosing
this bundle to be part of (some large power) of projective vector bundle
data these local maps may be made into global smooth homomorphism into the
image bundle
\begin{equation}
f\co E_\tau^N\longrightarrow\Cinf(M/X;E^-)
\text{ such that }P+f\text{ is surjective.}
\label{mms2001.366}\end{equation}
This necessarily stabilizes the null bundle to projective
vector bundle data with respect to the trivialization and we set
\begin{equation}
\index_a(P)=\left[\ker(P+f)-E_\tau^N\right]\in K^0(X,\cA).
\label{mms2001.367}\end{equation}
As in the untwisted case this class can be seen to be independent of the
precise stabilization used and to be homotopy invariant. In fact adding a
further stabilizing bundle is easily seen to leave the index unchanged and
stabilizing the family with the additional parameter of a homotopy shows
the homotopy invariance.

\begin{proposition}\label{mms2001.322} For a fibration \eqref{mms2001.320},
Azumaya bundle $\mathcal{A}$ over $X$ and projective bundles $E^\pm,$
there is a quantization of a given ismorphism $b$ of the lifts of these
bundles to $S^*(X/M)$ for which the null spaces, and hence also the null
spaces of the adjoint family, define a projective bundle over the base so
that the difference class $\index_a(b)\in K^0(X,\mathcal{A})$ depends only
on the class of $b$ in $K^0(T(X/M),\rho ^*\phi ^*\mathcal{A})$ and so
defines the analytic index homomorphism
\begin{equation}
\index_a\co  K (T(M/X),\rho^* \phi^*\mathcal{A}) \longrightarrow K
(X;\mathcal{A}).
\label{mms2001.323}\end{equation}
\end{proposition}

\begin{proof} The stabilization discussed above associates to an elliptic
    family with principal symbol $b$ an element of $K(X,\mathcal{A}).$ This
    class is independent of the stabilization used to define it and is
    similarly independent of the quantization chosen, since two such families
    differ by a family of compact operators. Clearly the element is unchanged
    if the symbol, or operator, is stabilized by the identity on some other
    primitve vector bundle defined over the lift of the same $\PU(N)$
    bundle. Furthermore the homotopy invariance of the index and the existence
    of a quantization map show that the element depends only on the homotopy
    class of the symbol. The additivity of the index under composition and the
    multiplicativity of the symbol map then shows that the resulting map
    \eqref{mms2001.323} is a homomorphism.
\end{proof}

\section{The topological index}\label{Topological}

In this section we define the topological index map for a fibration of compact
manifolds \eqref{mms2001.320}
\begin{equation}
\index_t\co \Kc(T(M/X),\rho ^*\phi ^*\mathcal{A})\longrightarrow
K^0(X,\mathcal{A})
\label{mms2001.325}\end{equation}
where $\rho \co T(X/M) \longrightarrow M$ is the projection and $\mathcal{A}$
is an Azumaya bundle over $X.$

We first recall some functorial properties of twisted $K$--theory. Let
$f\co Y\longrightarrow Z$ be a smooth map between compact manifolds.
Then the pullback map,
$$
f^!\co K(Z,\mathcal{A}) \longrightarrow K(Y, f^*\mathcal{A}),
$$
for any Azumaya bundle $\mathcal{A},$ is defined as follows. Let $V$
be finite dimensional projective vector bundle data over $Z,$ associated
with a trivialization of $\mathcal{A}.$ Then $f^* V $ is projective
vector bundle data over $Y$ associated to the lifted
trivialization and the resulting class in $K$--theory is independent of
choices, so defines $f^!.$ Alternatively, if $s$ is a section of the
twisted Fredholm bundle of $\mathcal{A}\otimes\mathcal{K}$ over $Z$, then
the pullback $f^*s$ is a section of the corresponding twisted Fredholm
bundle over $Y.$ The pull-back map may also be defined directly in terms of
the pull-back of projective bundles from the $\PU(N)$ bundle associated to
$\mathcal{A}$ over $Z$ to its pull-back over $Y.$

\begin{lemma}\label{Bott} For any Azumaya bundle there is a
canonical isomorphism
$$
j_! \co  K(X,\mathcal{A})\cong \Kc (X\times \R^{2N}, p_1^*\mathcal{A})
$$
determined by Bott periodicity.
\end{lemma}

\begin{proof} Recall that $\Kc (X\times \R^{2N}, p_1^*\mathcal{A}) \cong
K(C_0(X\times \R^{2N}, \E_{p_1^*\mathcal{A}})).$ Now $\E_{p_1^*\mathcal{A}}
\cong p_1^*\E_{\mathcal{A}},$ so that $C_0(X\times \R^{2N},
\E_{p_1^*\mathcal{A}}) \cong C(X,\E_{\mathcal{A}}) \widehat\otimes
C_0(\R^{2N}).$ Thus,
$$
\Kc (X\times \R^{2N}, p_1^*\mathcal{A}) \cong K(C(X,\E_{\mathcal{A}}) \otimes C_0(\R^{2N})).
$$
Together with Bott periodicity, $K(C(X,\E_{\mathcal{A}}) \otimes
C_0(\R^{2N})) \cong K(X,\mathcal{A}),$ this proves the lemma.
\end{proof}

If $\phi\co M\longrightarrow X$ is our basic fibre bundle of compact manifolds
we know that there is an embedding $f\co M\longrightarrow X \times \R^{N}$,
cf \cite{MR43:5554} \S 3. Then the \emph{fibrewise} differential is an
embedding $Df\co T(M/X)\longrightarrow X\times\R^{2N}$ with complex normal
bundle. In the untwisted case we have, via the Thom isomorphism, $Df_1 \co K_c
(T(M/X)) \longrightarrow K_c (X\times\R^{2N}).$ By anology with the case of
compact $\spinC$ manifolds, we call this the Gysin map.

We explain the extension to the twisted case. So again let $\mathcal A$ be
an Azumaya algebra over $X$ with $Y= T(M/X).$ Let $E$ be the (complex) normal
bundle to the imbedding $i\co M\longrightarrow X \times\RR^N$, and let $\mathcal
A_E$ be the lift of $\mathcal A$ to $E.$ Then 
\begin{equation}
\begin{aligned}
i_!\co  K_c(Y,\cA) &\longrightarrow \Kc (E,\cA_E),\\
\xi &\longmapsto (\pi^*\xi,\pi^*G) \otimes (\pi^* {S^+}, \pi^*{S^-}, c(v))
\end{aligned}
\label{Extra}\end{equation}
where $\xi=(\xi^+,\xi^-)$ is pair of projective vector bundle data over $X,$
associated to a local trivialization of $\cA,$ with $G\co \xi^+\longrightarrow
\xi^-$ a projective bundle map between them which is an isomorphism outside a
compact set and $(\pi^* {S^+}, \pi^*{S^-}, c(v))$ is the usual Thom class
of the complex vector bundle $E.$ On the the right hand side the
the graded pair of projective vector bundle data is 
\begin{equation*}
(\pi^*\xi^+\otimes\pi^*S^+\oplus\pi^*\xi^-\otimes\pi^*S^-,
\pi^*\xi^+\otimes\pi^*S^-\oplus\pi^*\xi^-\otimes\pi^*S^+)
\label{Extra2}\end{equation*}
with map between them being
\begin{equation*}
\left[\begin{matrix}
G&c(v)\\ c(v)&G^{-1}
\end{matrix}\right],\ v\in E.
\end{equation*}
This is an isomorphism outside a compact subset of $E$ and defines
a class in $\Kc(E,\cA_E)$ which is independent of choices. The {Thom isomorphism}
in this context, cf \cite{Donovan-Karoubi1}, asserts that $\;i_!\;$ is
an isomorphism.

Now, $E$ is diffeomorphic to a tubular neighborhood $\U$ of the image of
$Y;$ let $\Phi\co  \U \longrightarrow E$ denote this diffeomorphism. By the
Thom isomorphism above 
$$
i_!\co K_c(Y,i^*\mathcal{A})\longrightarrow \K_c (E,\mathcal{A}_E)\cong
\K_c (\U,\mathcal{A}'),
$$
where $\mathcal{A}'= \Phi^*(\mathcal{A}_E).$  Using Lemma~\ref{mms2001.327}, the
inclusion of the open set $U$ in $X \times R^{2N}$ induces a map
$\Kc (\U,\mathcal{A}')\longrightarrow\Kc(X\times R^{2N},
\pi_{1}^{*}\mathcal{A})$ where $\pi_1 \co  X\times \R^{2N}\longrightarrow X$
is the projection. The composition of these maps defines the Gysin
map. In particular we get the Gysin map in twisted $K$--theory,
$$
Df_! \co \Kc(T(M/X), \rho^*\pi^*\mathcal{A}) \longrightarrow\Kc(X\times
\R^{2N},
\pi_1^*\mathcal{A})
$$
where $\rho\co T(M/X)\longrightarrow X$ is the projection map. Since
$\pi = \pi_1 \circ f$ it follows that $Df^*\pi_1^*\mathcal{A} =
\rho^*\pi^*\mathcal{A}$. Now define the \emph{topological index},
\eqref{mms2001.325} as the map
$$
\index_t=j_!^{-1}\circ Df_! \co  \Kc(T(M/X), \rho^*\phi^*\mathcal{A})
\longrightarrow K(X,\mathcal{A}),
$$
where we apply Lemma~\ref{Bott} to see that the inverse $j_!^{-1}$ exists.

\section{Proof of the index theorem in twisted $K$--theory}\label{Equality}

We follow the axiomatic approach of Atiyah--Singer to prove that the
analytic index and the topological index coincide.

\begin{definition*}\label{mms2001.335}
An \emph{index map} is a homomorphism
\begin{equation}
\index\co \Kc(T(M/X),\rho^*\pi^*\mathcal{A})\longrightarrow K(X,\mathcal{A}),
\label{mms2001.356}\end{equation}
satisfying the following:

\begin{enumerate}
\item (Functorial axiom)\qua If $M$ and $M'$ are two fibre bundles with compact
fibres over $X$ and $f\co  M\longrightarrow M'$ is a diffeomorphism which
commutes with the projection maps $\phi\co  M\longrightarrow X$ and
$\phi'\co M'\longrightarrow X$ then the diagram
\begin{equation}
\xymatrix{\Kc(T(M/X),\rho^*\phi^*\mathcal{A})\ar[dr]_{\index}&&
\Kc(T(M'/X),\rho^*{(\phi')}^*\mathcal{A})\ar[dl]^{\index}\ar[ll]^{f^!}\\
&K(X,\mathcal{A})&}
\label{mms2001.331}\end{equation}
is commutative.

\item (Excision axiom)\qua  Let $\phi\co M\longrightarrow X$ and
$\phi'\co M'\longrightarrow X$ be two fibre bundles of compact manifolds, and
let $\alpha :\U \subset M$ and $\alpha ':\U' \subset M'$ be two open
sets with a
diffeomorphism $g\co \U\cong\U'$ satisfying $\phi'\circ g=\phi$ used to
identify them, then the diagram
\begin{equation}
\xymatrix @=4pc @ur {
       \Kc (T(\U/\phi(\U)), \rho^*\alpha ^*\phi^*\mathcal{A})
\ar[r]^{\alpha _!}\ar[d]_{(\alpha')_!} &
\Kc(T(M/X), \rho^*\phi^*\mathcal{A}) \ar[d]^{\index}  \\
\Kc(T(M'/X), \rho^*{(\phi')}^*\mathcal{A})\ar[r]_{\index'} &
K(X,\mathcal{A})}
\label{mms2001.332}\end{equation}
is commutative.

\item (Multiplicativity axiom)\qua Let $V$ be a real vector space and suppose
that $i\co M\longrightarrow X\times V$ is an embedding which intertwines the
projection maps $\phi\co  M\longrightarrow X$ and $\pi_1\co  X\times
V\longrightarrow X$, ie $\pi_1\circ i=\phi$; the fibrewise differential
$i_*\co T(M/X)\longrightarrow X\times TV$ also intertwines the
projections. The one-point compactification $S(V\oplus \R)$ of $V$ is a
sphere with the canonical inclusion $e\co TV\longrightarrow TS(V\oplus \R)$
inducing the inclusion $e'=\Id\times e\co  X\times TV \longrightarrow
X\times TS(V\oplus \R).$ Then the
diagram
\begin{equation}
\xymatrix @=4pc  {
& \Kc(X \times TV,\pi_1^*\mathcal{A})\ar^{(e')_!}[d]\\
\Kc(T(M/X), \rho^*\phi^*\mathcal{A}) \ar[dr]_{\index}\ar^{(i_*)_!}[ur] &
\Kc(X\times TS(V\oplus\R), \pi_1^*\mathcal{A})\ar^{\index}[d]\\
& K(X,\mathcal{A})}
\label{mms2001.333}\end{equation}
commutes.

\item (Normalization axiom)\qua 
If the fibre bundle of compact manifolds $\phi\co  M$ $\longrightarrow X$ has
single point fibres, then the index map
\begin{equation}
\index\co  \Kc(T(M/X),\rho^*\pi^*\mathcal{A})=K(X,\mathcal{A})\longrightarrow
K(X,\mathcal{A})
\label{mms2001.334}\end{equation}
is the identity map.
\end{enumerate}
\end{definition*}

The next theorem asserts in particular that such an index map does exist.

\begin{thm}\label{topindex}
The topological index, $\index_t,$ is an index map.
\end{thm}

\begin{proof} We proceed to check the axioms above in turn.

If $f\co  M\longrightarrow M'$ is a diffeomorphism as in the statement
of the functoriality axiom, let $i\co  M' \longrightarrow X\times V$ be an
embedding commuting with the projections, where $V$ is a finite
dimensional vector space. Then $ i \circ f \co  M \longrightarrow
X\times V$ is also such
an embedding. Using these maps, we may identify the topological index maps as
$\index_t' = j_!^{-1} \circ (i_*)_! $ and $\index_t
= j_!^{-1} \circ ((i\circ f)_*)_! = j_!^{-1} \circ (i_*)_! \circ (f_*)_! =
\index'_t \circ (f_*)_!$, where $j \co  X \hookrightarrow X\times V$ is the zero
section embedding. Then the diagram
$$
\xymatrix @=7pc @ur {
                    \Kc(T(M/X),
\rho^*\phi^*\mathcal{A})\ar[d]_{\index_t}\ar[dr]^{((f_*)_!)^{-1}}
                      \\
                    K(X,\mathcal{A}) \ar[r]_{\index_t} & \Kc(T(M'/X),
\rho^*{\phi'}^*\mathcal{A})}
$$
commutes. Since $f$ is a diffeomorphism, $(f_*)_!)^{-1} = {f}^!$, which
establishes that $\index_t$ is functorial.

Next consider the excision axiom and let $i'\co  M'\longrightarrow X\times V$
be an embedding, and $j\co  X\hookrightarrow X \times V$ be the zero section
embedding. Then $\index'_t = j_!^{-1} \circ (i_*)_!$, so that the relevant
map in the lower part of the diagram \eqref{mms2001.332} is $ j_!^{-1}
\circ (i_*)_! \circ (\alpha'_*)_!  = j_!^{-1} \circ ((i \circ
\alpha')_*)_!.$ But this is merely the topological index $\Kc
(T(\U/\pi(\U)), \rho^*\alpha^*\pi^*\mathcal{A}) \longrightarrow K(X,
\mathcal{A}).$ That it agrees with the relevant map in the upper part of
the diagram, follows from the fact that the topological index is well
defined and independent of the choice of embedding. Thus the topological
index satisfies the excision property.

The multiplicativity property for the topological index follows from its
independence of the choice of embedding, since \eqref{mms2001.333} amounts
to the definition of the topological index for the given embedding. The
independence of the choice of embedding is established briefly as
follows. Let $i_k\co  M\longrightarrow X\times V_k$, $k=0, 1$ be two
embeddings, and $j_k\co  X\hookrightarrow X \times V_k$ be the zero section
embedding, $k=0, 1$. Consider a linear homotopy $I_t\co  M\longrightarrow
X\times V_0\oplus V_1$, $t\in [0,1]$ defined as $I_t(m) = (i_0(m), t
i_1(m))$, and the zero section embedding $J \co  M\longrightarrow X\times
V_0\oplus V_1$ defined as $J(m) = (j_0(m), j_1(m))$.  By the homotopy
invariance of the induced map in $K$--theory, it follows that $J_!^{-1}
\circ ({I_t}_*)_!$ is independent of $t$. Using the functorial 
property 
of the topological index, one deduces that
${j_0}_!^{-1} \circ ({i_0}_*)_!$ agrees with $J_!^{-1} \circ
({I_1}_*)_!$. Now let  $\tilde I_t \co  M\longrightarrow
X\times V_0\oplus V_1$, $t\in [0,1]$ be defined as $\tilde I_t(m) = (t i_0(m),
i_1(m))$. The argument above establishes that ${j_1}_!^{-1} \circ
({i_1}_*)_!$ also agrees with $J_!^{-1} \circ (\tilde{ I_1}_*)_!$ = 
$J_!^{-1} \circ ({I_1}_*)_!$. Therefore ${j_0}_!^{-1} \circ 
({i_0}_*)_!
= {j_1}_!^{-1} \circ ({i_1}_*)_!$ as claimed.

For the normalization axiom, note that in case $M=X,$ if
$i\co X\longrightarrow X\times V$ is an embedding which commutes with the
projection maps $\phi\co  X\longrightarrow X$ and $\pi_1\co  X\times
V\longrightarrow X$ then $i$ is necessarily the trivial embedding
$\Id\times g$ with $g\co X\longrightarrow V$ constant. Then $\index_t =
j_!^{-1} \circ (i_*)_! = \Id$ since $i_* = \Id \times 0,$ which shows that
the topological index is normalized.
\end{proof}

\begin{thm}\label{indexunique} There is a unique index map.
\end{thm}

\proof We have already shown that $\index_t$ is an index map. Thus
it suffices to consider a general index map as in \eqref{mms2001.356} and
to show that $\index=\index_t.$

Suppose that $i\co  M\longrightarrow X\times V$ is an embedding
which intertwines the projection maps $\phi\co  M\longrightarrow X$ and
$\pi_1\co  X\times V\longrightarrow X.$ Then, together with the notation of
the Multiplicativity Axiom, set $i^+= e'\circ i_*\co T(M/X)\longrightarrow
X\times TS(V\oplus \R).$ Let $0\in TV$ be the origin and
$j\co \{0\}\longrightarrow TV$ be the inclusion, inducing the inclusion
$j'=\Id\times j\co X\times\{0\}\longrightarrow X\times TV$ and denote the
composite inclusion $j^+=e'\circ j'\co X\times\{0\}\longrightarrow X\times
TS(V\oplus\R).$ Then consider the diagram:
$$
\xymatrix @=2.9pc{& \Kc(T(M/X), \rho^*\phi^*\mathcal{A})
\ar[dl]_{(i_*)_!}\ar[d]^{i^+}\ar[dr]^{\index} & &\\ \Kc(X\times TV,
\pi_1^*\mathcal{A}) \ar^{e_!\qquad}[r] &\Kc(X\times TS(V\oplus \R),\pi_1^*\mathcal{A}
)\ar[r]\ar[r]^{\qquad\quad \index^S} & K(X, \mathcal{A}) \\ & K(X,
\mathcal{A})\ar[ur]_{\index} \ar[ul]^{j_!}\ar[u]_{j^+}  & }
$$
The left side of this diagram commutes by the excision property and the
right side by the multiplicative property. By the normalization property,
the composite map $\index\circ j^+$ is the identity mapping, so
\begin{align*}
\index =& \index^{S}\circ i^+_!= \index^{S}\circ e_!\circ i_!=
\index^{S}\circ j^+_!\circ j^{-1}_!\circ i_!\\
=& \index'\circ j^{-1}_!\circ i_!= j^{-1}_!\circ i_!=\index_t.\tag*{\qed}
\end{align*}

The following theorem completes the proof of the index theorem in twisted
$K$--theory.

\begin{thm}\label{analyticindex}The analytic index $\index_a$ is an index map.
\end{thm}

\begin{proof} Again we consider the axioms for an index map in order.

The invariance of the algebra of pseudodifferential operators under
diffeomorphism, and the naturality in this sense of the symbol map, show
that under the hypotheses of the functoriality axiom, there is an isomorphism
of short exact sequences \eqref{mms2001.310}:
{\small
\begin{equation}
\begin{CD}
\Psi^{-1}(M'/X, \E) @>>> \Psi^{0}(M'/X, \E) @>>>
S^{[0]}(T(M'/X),\rho^*\End(\E))
\\
                @VV{f^*}V        @VV{f^*}V             @VV{f^*}V&
\\ \Psi^{-1}(M/X, \E) @>>> \Psi^{0}(M/X, \E) @>>>
S^{[0]}(T(M/X),\rho^*\End(\E)).
\end{CD}
\end{equation}}
Since the analytic index is by definition the boundary map in the associated
6--term exact sequence in $K$--theory, we see that $\index_a(f^*[p]) =
f^!\index_a([p])$, for all $[p] \in \Kc(T(M/X), \rho^*\pi^*\mathcal{A}).$ This
establishes the functoriality of $\index_a.$

For the excision axiom, observe that any element in $\Kc(T(\U/\pi(\U)),
\rho^*i^*\pi^*\mathcal{A})$ may be represented by a pair of projective
vector bundle data over $\U$ and a symbol $q \in S^{[0]}(T(\U/\pi(\U)$ with
the property that $q$ is equal to the identity homomorphism outside a compact
set in $\U.$ Complementing the second bundle with respect to vector bundle
over $M,$ using the discussion following Lemma~\ref{mms2001.352}, we may
extend both sets of projective vector data to the whole of $M,$ to be equal
outside $\U.$ This also extends $q$ to an element $p\in S^{[0]}(T(M/X)$ by
trivial extension. The exactness in \eqref{mms2001.310} shows that there is
a projective family of elliptic pseudodifferential operators $P$ of order
zero with symbol equal to $p,$ by use of a partition of unity we may take
it to be equal to the identity outside $\U$ in $M.$ Similarly, $q$ also
defines an element $p'\in S^{[0]}(T(M'/X)$ and, from the corresponding
exact sequence, there is a projective family of elliptic pseudodifferential
operators $P'$ equal to the identity outside $\U$ in $M';$ we may further
arrange that $P=P'$ in $\U.$ We can construct parametrices $Q$ of $P$ and
$Q'$ of $P'$ such that $Q$ is equal to the identity outside $\U$ in $M$ and
$Q'$ is equal to the identity outside $\U$ in $M'$ and $Q=Q'$ in $\U.$ By
the explicit formula for the analytic index in terms of the projective
family of elliptic pseudodifferential operators and its parametrix, see \S
3, it follows that the diagram \eqref{mms2001.332} commutes, that is, the
analytic index satisfies the excision property.

Under the hypotheses of the multiplicative axiom we need to show for a
class $[p]\in \Kc(T(M/X), \rho^*\phi^*\mathcal{A}),$ represented by a
symbol $p\in S^0(T(M/X),\E),$ that $\index_a([p]) = \index_a(h_![p]),$
where $h\co  T(M/X) \longrightarrow X \times TS(V\oplus \R)$ is the embedding
that is obtained as the composition $h = e\circ Di$, and $h_!$ is the Gysin
map. This is done by first embedding $M$ as the zero section of the
compactification of its normal bundle to a sphere bundle. In this case
one may argue as in \cite{MR43:5554}, where a family of operators $B$ is
constructed on the sphere $S^n=S(\R^n\times\R)$ to be $O(n)$ invariant,
surjective and have symbol equal to the Thom class. Then $B$ can be
extended naturally to act on the fibres of the sphere bundle. Having
stabilized $P,$ to a projective family with the given symbol $p,$ (and a
finite rank term $f)$ we may lift it, as described in \cite{MR43:5554}, to
be an operator acting on the lift of $E^\pm$ to the sphere bundle and
reducing to $P$ on fibre-constant sections. As in \cite{MR43:5554} the
tensor product of the lifted operator $P$ and $B$ then acts as a Fredholm
family
\begin{equation}
\begin{pmatrix}P&B\\B^*&P^*\end{pmatrix}.
\label{mms2001.368}\end{equation}
between the bundles $E^+\otimes G^+\oplus E^-\otimes
G^-$ and $E^-\otimes G^+\oplus E^+\otimes G^-.$ Since $P$ and $B$ commute
it follows as in the untwisted case that the null space of this surjective
operator is isomorphic to the null space of $P.$ Thus has represents the
same index class which proves the desired multiplicativity in this case.

The general case now follows by using the excision property, so the analytic
index satisfies the multiplicative property.

The normalization axiom holds by definition; it is important that this is
consistent with the proof of the axioms above.
\end{proof}

The equality of the topological and analytic indexes is now an immediate
consequence of Theorems \ref{topindex}, \ref{indexunique} and
\ref{analyticindex}:

\begin{thm}[The index theorem in $K$--theory] \label{K-indexthm} Let
$\phi\co M\longrightarrow X$ be a fibre bundle of compact manifolds, let
$\mathcal{A}$ be an Azumaya bundle over $X$ and $P$ be a projective family
of elliptic pseudodifferential operators acting between two sets of
projective vector bundle data associated to a local trivialization of
$\mathcal{A}$ and with symbol having class $p\in\Kc(T(M/X);\rho ^*\phi
^*\mathcal{A}),$ where $\rho\co T(M/X)\longrightarrow M$ is the projection
then \begin{equation}\label{Kind}
\index_a(P) = \index_t(p) \in K(X,\mathcal{A}).
\end{equation}
\end{thm}

\section{The Chern character of the index bundle}\label{Chern}

In this section, we compute the Chern character of the index bundle and
obtain the cohomological form of the index theorem for projective families
of elliptic pseudodifferential operators.  In the process, not
surprisingly, the torsion information from the Azumaya bundle is lost. We
begin with the the basic properties of the Chern character.

The Chern character of projective vector bundles, defined in
\S\ref{SS.Chern}, gives a homomorphism
\begin{equation}\label{Chern-char}
\Ch_\A \co K^0(X, \mathcal{A})\longrightarrow H^{even}(X, \bbQ).
\end{equation}
It satisfies the following properties.

\begin{enumerate}
\item The Chern character is \emph{functorial} under smooth maps
in the sense that if $f\co  Y \longrightarrow X$ is a smooth map between
compact manifolds, then the following diagram commutes:
\begin{equation} \label{natural2}
\begin{CD}
K^0(X, \mathcal{A}) @>f^!>> K^0(Y, f^*\mathcal{A}) \\
                @VV{\Ch_\A}V
@VV{\Ch_{f^*\A}}V&
\\H^{even}(X, \Q	) @>f^*>>   H^{even}(Y, \Q).
\end{CD}
\end{equation}

\item The Chern character respects the module structure, of $K^0(X,
\mathcal{A})$ over $K^0(X),$ in the sense that the following diagram
commutes:
\begin{equation} \label{module}
\begin{CD}
K^0(X) \times K^0(X, \mathcal{A}) @>>> K^0(X, \mathcal{A}) \\
                @VV{\Ch \times \Ch_\A}V
@VV{\Ch_{\A}}V&
\\H^{even}(X, \Q	) \times H^{even}(X, \Q	) @>>>   H^{even}(X, \Q	)
\end{CD}
\end{equation}
where the top horizontal arrow is the action of $K^0(X)$ on $K(X,\mathcal{A})$
given by tensor product and the bottom horizontal arrow is given by
the cup product.
\end{enumerate}

\begin{thm}[The cohomological formula of the index theorem]
\label{Chernindex}
For a fibration \eqref{mms2001.320} of compact
manifolds and a
projective family of elliptic pseudo-\break 
differential operators $P$ with symbol
class $p \in \Kc(T(M/X), \rho^*\phi^*\mathcal{A})$, where\break $\rho \co 
T(M/X)\longrightarrow M$ is the projection, then
\begin{equation}\label{cohind}
\Ch_\A (\index_a P) = (-1)^n
\phi_*\rho_*\left\{\rho^*{\rm Td}(T(M/X)\otimes
\C)\cup\Ch_{\rho^*\phi^*\A}(p)
\right\}
\end{equation}
where the Chern character is denoted $\Ch_\A \co  K(X, \A) \to H^\bullet (M)$ and
$\Ch_{\rho^*\phi^*\A}\co$  $\Kc(T(M/X), \rho^*\phi^*\mathcal{A})
\to H^\bullet_c(T(M/X)),$
$n$ is the dimension of the fibres of $\phi\circ\rho$,
$\phi_*\rho_*\co  H^\bullet_c(T(M/X))\longrightarrow H^{\bullet-n}(X)$ is
integration
along the fibre.
\end{thm}

This theorem follows rather routinely from the index theorem in $K$--theory,
Theorem~\ref{K-indexthm}. The key step to getting the formula
\eqref{cohind} is the analog of the Riemann--Roch formula in the context of
twisted $K$--theory, which we now discuss.

Let $\pi\co  E\longrightarrow X$ be a $\spinC$ vector bundle over $X$,
$i\co  X\longrightarrow E$ the zero section embedding,
$F$ be a complex projective vector bundle over $X$ that is associated
to the Azumaya bundle $\A$ on $X$. Then we compute,
$$
\begin{array}{lcl}
\Ch_{\pi^*\A}(i_! F) &=& \Ch_{\pi^*\A}(i_! 1 \otimes \pi^*F)\\[+7pt]
		 &=& \Ch(i_! 1) \cup \Ch_{\pi^*\A}(\pi^* F),
\end{array}
$$ where we have used the fact that the Chern character
respects the $K^0(X)$--module structure.
The standard Riemann--Roch formula asserts that
$$
\Ch(i_! 1) = i_*{\rm Td}(E)^{-1}.
$$
Therefore we obtain the following Riemann--Roch formula for twisted
$K$--theory,
\begin{equation}\label{RR}
\Ch_{\pi^*\A}(i_! F) = i_*\left\{{\rm Td}(E)^{-1} \cup \Ch_\A( F)
\right\}.
\end{equation}

The index theorem in $K$--theory in \S 5 shows in particular that
$$
\Ch_\A (\index_a P) = \Ch_\A (\index_t p).
$$
Now $\index_t p = j_!^{-1} \circ (Di)_! $ where
$i  \co  M \hookrightarrow X\times V$ is an
embedding that commutes with the projections
$\phi \co  M\longrightarrow X$ and $\pi_1 \co  X\times V \longrightarrow X$, and
$j \co  X \hookrightarrow X\times V$ is the zero section
embedding.
Therefore
$$
\Ch_\A (\index_t p) = \Ch_\A (j_!^{-1} \circ (Di)_! p)
$$
By the Riemann--Roch formula for twisted $K$--theory \eqref{RR},
$$
\Ch_{\pi_1^* \A}(j_! F) =  j_*\Ch_\A(F)
$$
since $\pi_1\co  X\times V \longrightarrow X$ is a trivial bundle.
Since
${\pi_1}_* j_* 1 = (-1)^n$, it follows that
for $\xi \in \Kc(X \times V, \pi_1^*\A) $, one has
$$
\Ch_\A(j_!^{-1} \xi) = (-1)^n{\pi_1}_*\Ch_{\pi_1^*\A}( \xi)
$$
Therefore
\begin{equation}\label{a1}
\Ch_\A (j_!^{-1} \circ (Di)_! p)
=  (-1)^n {\pi_1}_*\Ch_{\pi_1^*\A}((Di)_! p)
\end{equation}
By the Riemann--Roch formula for twisted
$K$--theory \eqref{RR},
\begin{equation}\label{a2}
\Ch_{\pi_1^*\A}((Di)_!  p) = (Di)_*\left\{\rho^*{\rm Td}(N)^{-1}
\cup\Ch_{\rho^*\phi^*\A}( p)\right\}
\end{equation}
where $N$ is the complexified normal bundle to the embedding
$Di \co  T(M/X) \longrightarrow X \times TV$, that is
$N =  X \times TV/Di(T(M/X))\otimes \C$.
Therefore ${\rm Td}(N)^{-1} = {\rm Td}(T(M/X)\otimes \C)$ and
          \eqref{a2} becomes
$$
\Ch_{\pi_1^*\A}((Di)_! p) =(Di)_*\left\{\rho^*{\rm Td}(T(M/X)\otimes \C)
\cup\Ch_{\rho^*\phi^*\A}( p)\right\}.
$$
Therefore \eqref{a1} becomes
\begin{equation}
\begin{array}{lcl}
\Ch_\A (j_!^{-1} \circ (Di)_! p) &= &
(-1)^n {\pi_1}_*(Di)_*\left\{\rho^*{\rm Td}(T(M/X)\otimes \C)  \cup
\Ch_{\rho^*\phi^*\A}( p)\right\}\\[+7pt]
& = & (-1)^n \phi_*\rho_*\left\{\rho^*{\rm Td}(T(M/X)\otimes \C)  \cup
\Ch_{\rho^*\phi^*\A}( p)\right\}
\end{array}
\end{equation}
since $\phi_* \rho_*= {\pi_1}_* (Di)_*$. Therefore
\begin{equation}\label{cohtopindex}
\Ch_\A (\index_t p) = (-1)^n \phi_*\rho_*\left\{\rho^*{\rm
Td}(T(M/X))
\cup
\Ch_{\rho^*\phi^*\A}( p)\right\},
\end{equation}
proving Theorem~\ref{Chernindex}.

\section{Determinant line bundle of the index bundle}\label{Determinant}

In this section, we define the determinant line bundle of the index bundle
of a projective family of elliptic pseudodifferential operators, and
compute its Chern class.

We begin with the general construction of the determinant line
bundle of a projective vector bundle over $X$.
Let $\A$ be an Azumaya
bundle over $X$ and
    $P\stackrel{\pi}{\to} X$ be the principal ${\rm
PU}(n)$ bundle associated to $\A$.
Let $E\to P$ be a projective
vector bundle
over $X$, associated to $\A$, cf section 1.3.
Recall that $E$
satisfies in addition the condition
\begin{equation} \label{bg}
\phi^*E \cong \pi_P^*E \otimes \pi_{{\rm PU}(n)}^*L
\end{equation}
where $\phi \co  {\rm PU}(n)\times P \to P$ is the action, $\pi_P \co 
{\rm PU}(n)\times P \to P$ is the projection onto the second factor,
      $\pi_{{\rm PU}(n)} \co   {\rm PU}(n)\times P \to {\rm PU}(n)$
is the projection onto the first factor, $L \to {\rm PU}(n)$ is the
(determinant) line bundle associated to the principal ${\mathbb Z}_n$ bundle
${\mathbb Z}_n\to {\rm SU}(n)\to {\rm PU}(n)$ as in section 1.3.

Then we observe that
$$
\phi^*\Lambda^n (E ) \cong \pi_P^*\Lambda^n (E)
$$
and therefore $\Lambda^n (E ) = \pi^*(F)$ for some line bundle
$F\to X.$ Define $\det(E) = F$ to be the determinant line bundle of
the projective vector bundle $E$.

     In particular, this gives a homomorphism
$$
\det \co  K^0(X, \mathcal{A}) \longrightarrow \pi_0({\rm Pic}(X))
$$
where ${\rm Pic}(X)$ denotes the Picard variety of $X,$ and the
components of the Picard variety $\pi_0({\rm Pic}(X))$ consist of the
isomorphism classes of complex line bundles over $X.$

\begin{thm}[Chern class of the determinant line bundle of the index
bundle] For a fibration \eqref{mms2001.320} of compact manifolds and a
projective family of elliptic pseudodifferential operators $P$ with symbol
class $p \in \Kc(T(M/X), \rho^*\phi^*\mathcal{A})$,  where $\rho \co 
T(M/X)\longrightarrow M$ is the projection,
\begin{equation}\label{cohind2}
{\rm c}_1(\det (\index_a P) )= \{(-1)^n
\phi_*\rho_*\left\{\rho^*{\rm Td}(T(M/X)\otimes
\C)\cup\Ch_{\rho^*\phi^*\A}(p)
\right\}\}^{[2]}
\end{equation}
where $\Ch_{\rho^*\phi^*\A} \co  \Kc(T(M/X), \rho^*\phi^*\A) \longrightarrow
H^\bullet_c(T(M/X))$ is the Chern character, $c_1$ is the first Chern
class, $N$ is the dimension of the fibres of $\phi\circ\rho,$
$\phi_*\rho_*$ is integration along  the fibres mapping
$H^\bullet_c(T(M/X))$ to $H^{\bullet-n}(X)$ is and $\{\cdot\}^{[2]}$
denotes the degree 2 component.
\end{thm}

\begin{proof}
The proof of the theorem follows from Theorem~\ref{Chernindex} and the
second part of Lemma~\ref{chern}.
\end{proof}

\section{Projective families of Dirac operators}\label{Dirac}

Let $\Cl(M/X)$ denote the bundle of Clifford algebras on the fibres of
$\phi$ for some family of fibre metrics.  A fiberwise Clifford module on a
bundle $E$ over the total space of a fibration is a smooth action of
$\Cl(M/X)$ on the bundle. That is to say it is an algebra homomorphism
\begin{equation}
\Cl(M/X)\longrightarrow \End(E).
\label{mms2001.369}\end{equation}
Since the endomorphism bundle of a projective bundle over $E,$ with respect
to an Azumaya bundle $\cA,$ is a vector bundle, this definition can be
taken directly over to the projective case. Similarly, the condition that
the Clifford module structure be hermitian can be taken over as the
condition that \eqref{mms2001.369} be $*$--preserving. The condition that a
unitary connection on $E$ be a Clifford connection is then the usual
distribution condition, for vertical vector fields,
\begin{equation}
\nabla_X(\cli(\alpha)e_a)=\cli(\nabla_X\alpha )e_a+\cli(\alpha )(\nabla_Xe_a)
\label{mms2001.370}\end{equation}
in terms of the Levi--Civita connection on the fibre Clifford bundle and for
any sections $e_a$ of the bundle trivialization of $E$ with respect to a
full trivialization of $\cA.$

The Dirac operator associated to such a unitary Clifford connection on a
hermitian projective Clifford module is then given by the usual formula
over the open sets $U_a$ of a given full trivialization of $\cA$ over the
base:
\begin{equation}
\eth_ae_a=\widetilde\cli(\nabla e_a)
\label{mms2001.371}\end{equation}
where $\widetilde\cli$ is the contraction given by the Clifford action from
$T^*(M/X)\otimes E_a$ to $E_a.$ The invariance properties of the connection
and Clifford action show that the $\eth_a$ form a projective family of
differential operators on the projective bundle $E.$

As in the untwisted case, if $\eth$ is a Dirac operator in this sense,
acting on a vector bundle $F$ over $M$ and $E$ is a projective vector
bundle over $M,$ relative to some Azumaya algebra $\cA,$ then we may twist
$\eth$ by choosing a unitary connection $\nabla^E$ on $E$ and extending the
Clifford module trivially from $F$ to $F\otimes E$ (to act as multiples of
the identity on $E)$ and taking the tensor product connection on $F\otimes
E.$ The resulting Dirac operator is then a projective family as described
above.

In particular, if the bundle $T(M/X)$ is spin and we consider the
family of  Dirac operators along the fibres of $\phi$ twisted by the
projective vector bundle $E$ over $M,$ we deduce from Theorem \ref{K-indexthm}
that
$$
\index_a(\eth)= \phi_!(E)  \in
K(X, \A)
$$
where $\phi_!\co   K(M, \phi^*\A) \to K(X, \A)$ is defined as
$\phi_! =
j_!^{-1}\circ f_! $ in the notation of Section \ref{Topological}.
By
Theorem \ref{Chernindex} and standard
manipulations of
characteristic classes one has
$$
\Ch_\A(\index_a(\eth))
= \phi_*(\widehat{A}(M/X)\Ch_{\phi^*\A}(E)) \in H^\bullet(X).
$$
where $\phi_*\co   H^\bullet(M) \to H^{\bullet-\ell}(X)$ denotes
integration along the $\ell$--dimensional fibres. A similar formula holds
more generally, in case $T(M/X)$ is a $\spinC$ bundle, with an extra factor
of $\exp(c_1)$ arising from the twisting curvature, see
\cite{Berline-Getzler-Vergne1}.

\end{document}

%% file: gtspec.tex
\def\ifplaintex{\expandafter\ifx\csname documentclass\endcsname\relax}

\def\ifplaintex{\expandafter\ifx\csname documentclass\endcsname\relax}

\ifplaintex 
\else
\fi

\expandafter\ifx\csname epsfbox\endcsname\relax\input epsf\fi

\def\gt{{\mathsurround=0pt\it $\cal G\mskip-2mu$eometry \&\ 
$\cal T\!\!$opology}}        

\def\gtp{{\mathsurround=0pt\it $\cal G\mskip-2mu$eometry \&\ 
$\cal T\!\!$opology $\cal P\!$ublications}}  

\def\lognumber#1{\def\thelognumber{#1}}
\def\volumenumber#1{\def\thevolumenumber{#1}}
\def\papernumber#1{\def\thepapernumber{#1}}
\def\volumeyear#1{\def\thevolumeyear{#1}}

\def\pagenumbers#1#2{\def\startpage{#1}\def\finishpage{#2}}
\def\published#1{\def\publishdate{#1}}
\def\proposed#1{\def\theproposer{#1}}
\def\seconded#1{\def\theseconders{#1}}
\def\received#1{\def\receiveddate{#1}}

\def\accepted#1{\def\accepteddate{#1}}

\def\asciiaddress#1{\def\theasciiaddress{#1}}
\def\asciiemail#1{\def\theasciiemail{#1}}

\long\def\asciiabstract#1{\long\def\theasciiabstract{#1}}
\def\asciikeywords#1{\def\theasciikeywords{#1}}
\def\shortauthors#1{\def\theshortauthors{#1}}

\let\\\par\let\thelognumber\relax
\let\thevolumenumber\relax\let\thepapernumber\relax
\let\thevolumeyear\relax\let\thesamplenumber\relax\let\startpage\relax
\let\finishpage\relax\let\publishdate\relax\let\receiveddate\relax
\let\reviseddate\relax\let\accepteddate\relax\let\theasciititle\relax
\let\theasciiauthors\relax\let\theasciiaddress\relax
\let\theasciiabstract\relax\let\theasciikeywords\relax
\let\theasciiemail\relax\let\theshortauthors\relax\let\theshorttitle\relax

\long\def\maketitlep{   

\count0=\startpage

\gt\hfill      
\hbox to 77pt{\vbox to 0pt{\vglue -15pt\epsfbox{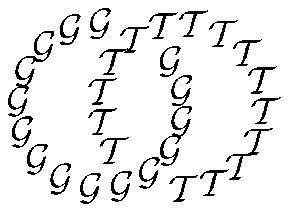}\vss}\hss}
\break
{\small\ifx\thesamplenumber\relax 
Volume \else Sample
\fi\thevolumenumber\ (\thevolumeyear)
\startpage--\finishpage\nl
Published: \publishdate}
\vglue 0.5truein plus 0.4fil minus 0.1truein

{\parskip=0pt\leftskip 0pt plus 1fil\def\\{\par\smallskip}{\ifplaintex\large
\else\Large\fi\bf\thetitle}\par\medskip}   

\vglue 0pt plus 0.1fil 

{\parskip=0pt\leftskip 0pt plus 1fil\def\\{\par}{\sc\theauthors}
\par\medskip}

\vglue 0pt plus 0.1fil 

{\small\parskip=0pt\let\newline\\
{\leftskip 0pt plus 1fil\def\\{\par}{\sl\theaddress}\par}
\expandafter\ifx\theemail\relax    
\relax\else\vglue 5pt plus 0.02fil minus 2pt\def\\{\stdspace{\rm 
and}\stdspace} 
\cl{Email:\stdspace\tt\theemail}\fi
\ifx\theurl\relax                  
\relax\else\vglue 5pt plus 0.02fil minus 2pt\def\\{\stdspace{\rm 
and}\stdspace}
\cl{URL:\stdspace\tt\theurl}\fi\par}

\vglue 7pt plus 0.3fil minus 3pt

{\bf Abstract}
\vglue 5pt plus 0.1fil minus 2pt

\theabstract

\vglue 7pt plus 0.3fil minus 3pt

{\bf AMS Classification numbers}\quad Primary:\quad \theprimaryclass

Secondary:\quad \thesecondaryclass

\vglue 5pt plus 0.3fil minus 2pt

{\bf Keywords:}\quad \thekeywords

\vglue 10pt plus 0.5fil minus 5pt

{\small  Proposed: \theproposer\hfill Received: \receiveddate\nl
Seconded: \theseconders\hfill 
\ifx\reviseddate\relax                         
Accepted: \accepteddate                        
\else
Revised: \reviseddate                          
\fi}
\eject
}       

\font\phead=cmsl9 scaled 950
\font=cmsl9 scaled 1050
\font\pnum=cmbx10 scaled 913
\font\lnum=cmbx10 
\font\pfoot=cmsl9 scaled 950
\font=cmsl9 scaled 1050
\ifplaintex
\headline{\vbox to 0pt{\vskip -4.5mm\line{\small\phead\ifnum
\count0=\startpage ISSN 1364-0380 (on line)
1465-3060 (printed) \hfill {\pnum\folio}\else\ifodd\count0\def\\{ }%
\ifx\theshorttitle\relax\thetitle\else\theshorttitle\fi\hfill{\pnum\folio}
\else\def\\{ and }{\pnum\folio}\hfill\ifx\theshortauthors\relax\theauthors
\else\theshortauthors\fi\fi\fi}\vss}}
\footline{\vbox to 0pt{\vglue 0mm\line{\small\pfoot\ifnum\count0=\startpage
\copyright\ \gtp\hfill\else
\gt, Volume \thevolumenumber\ (\thevolumeyear)\hfill\fi}\vss
}}
\else
\makeatletter
\def\@oddhead{{\small\lhead\ifnum\count0=\startpage ISSN 1364-0380 (on line)
1465-3060 (printed) \hfill {\lnum\number\count0}\else\ifodd\count0
\def\\{ }\ifx\theshorttitle\relax \thetitle \else\theshorttitle\fi\hfill
{\lnum\number\count0}\else\def\\{ and }{\lnum\number\count0}
\hfill\ifx\theshortauthors\relax 
\theauthors\else\theshortauthors\fi\fi\fi}}\def\@evenhead{\@oddhead}
\def\@oddfoot{\small\lfoot\ifnum\count0=\startpage\copyright\ \gtp\hfill\else
\gt, Volume \thevolumenumber\ (\thevolumeyear)\hfill\fi}
\def\@evenfoot{\@oddfoot}
\makeatother
\fi

\newwrite\gtoutfile
\long\gdef\makeheadfile{  
{\def\\{, }\def\s{ }
\immediate\openout\gtoutfile head.xxx
\immediate\write\gtoutfile{Proxy-for: \ifx\theasciiauthors\relax
\theauthors\else\theasciiauthors\fi\s<\ifx\theasciiemail\relax\theemail\else\theasciiemail\fi>}
\immediate\write\gtoutfile{\noexpand\\}
\immediate\write\gtoutfile{Authors: \ifx\theasciiauthors\relax
\theauthors\else\theasciiauthors\fi}
{\def\\{ }\immediate\write\gtoutfile{Title: \ifx\theasciititle\relax
\thetitle\else\theasciititle\fi}}
\immediate\write\gtoutfile{Subj-class: GT or SG or MG etc}
\immediate\write\gtoutfile{MSC-class: \theprimaryclass\ifx\thesecondaryclass\relax\else, \thesecondaryclass\fi}
\immediate\write\gtoutfile{Journal-ref: Geom. Topol. \thevolumenumber
(\thevolumeyear) \startpage-\finishpage}
\immediate\write\gtoutfile{Comments: Published by Geometry and Topology at}
\immediate\write\gtoutfile{\s\s http://www.maths.warwick.ac.uk/gt/GTVol\thevolumenumber/paper\thepapernumber.abs.html}
\immediate\write\gtoutfile{\noexpand\\}
\immediate\write\gtoutfile{}
\ifx\theasciiabstract\relax
\immediate\write\gtoutfile{\theabstract}\else
\immediate\write\gtoutfile{\theasciiabstract}\fi
\immediate\write\gtoutfile{}
\immediate\write\gtoutfile{\noexpand\\}
\immediate\write\gtoutfile{}
\immediate\closeout\gtoutfile}}  

\def\maketitlepage{\maketitlep\makeheadfile}
\let\maketitle\maketitlepage